\documentclass[oneside,12pt,letterpaper]{amsart}
\usepackage{amsmath,amsthm,amsfonts,amssymb}
\usepackage[alphabetic]{amsrefs}
\usepackage{hyperref}
\usepackage{subfigure}
\usepackage{eucal}
\usepackage{mathrsfs}
\usepackage{amssymb}

\usepackage[pdftex]{graphicx}

\usepackage{setspace}

\begin{document}

\theoremstyle{plain}
\newtheorem{thm}{Theorem}[section]
\newtheorem*{thm*}{Theorem}
\newtheorem{lem}[thm]{Lemma}
\newtheorem{prop}[thm]{Proposition}
\newtheorem{cor}[thm]{Corollary}

\theoremstyle{definition}
\newtheorem*{defn}{Definition}
\newtheorem*{exmp}{Example}

\theoremstyle{remark}
\newtheorem*{rem}{Remark}
\newtheorem{Rem}{Remark}
\newtheorem*{note}{Note}

\newcommand{\bignote}[1]{}

\newcommand{\bezout}{B\'{e}zout}

\newcommand{\N}{\mathbf{N}} 
\newcommand{\C}{\mathbf{C}} 
\newcommand{\Z}{\mathbf{Z}} 
\newcommand{\A}{\mathbf{A}}
\newcommand{\PR}{\mathbf{P}}
\newcommand{\into}{\hookrightarrow}
\newcommand{\Sym}{\text{Sym}}
\newcommand{\mult}{\text{mult}}
\newcommand{\PGL}{\text{PGL}}
\newcommand{\SL}{\text{SL}}
\newcommand{\Bl}{\text{Bl}}
\newcommand{\stab}{\text{stab}}
\newcommand{\codim}[1]{\text{codim}_{#1}~}
\newcommand{\Spec}{\text{Spec }}
\newcommand{\SPEC}{\text{\textbf{Spec} }}
\newcommand{\Proj}{\text{Proj}}
\newcommand{\Pic}{\text{Pic}}
\newcommand{\Span}{\mbox{Span }}
\newcommand{\fun}[3]{#1\!\!:#2\rightarrow #3}

\newcommand{\floor}[1]{\left\lfloor{#1}\right\rfloor\!}

\newcommand{\grd}[2]{\mathfrak{g}^{#1}_{#2}}

\newcommand{\isom}{\cong}

\newcommand{\scl}{\mathcal{G}}
\newcommand{\sclodd}{\mathcal{G}^{\text{odd}}}
\newcommand{\scleven}{\mathcal{G}^{\text{even}}}
\newcommand{\sclhyp}{\mathcal{G}^{\text{hyp}}}

\newcommand{\wpl}{\mathcal{C}}


\title{Subcanonical points on algebraic curves}
\author{Evan Merrill Bullock}


\maketitle

\begin{abstract}
If $C$ is a smooth, complete algebraic curve of genus $g\geq 2$ over the complex numbers, a point $p$ of $C$ is \emph{subcanonical} if $K_C \isom \mathcal{O}_C\big((2g-2)p\big)$.  We study the locus $\scl_g\subseteq \mathcal{M}_{g,1}$ of pointed curves $(C,p)$ where $p$ is a subcanonical point of $C$.  Subcanonical points are Weierstrass points, and we study their associated Weierstrass gap sequences.  In particular, we find the Weierstrass gap sequence at a general point of each component of $\scl_g$ and construct subcanonical points with other gap sequences as ramification points of certain cyclic covers and describe all possible gap sequences for $g\leq 6$.
\end{abstract}
\tableofcontents

\section{Introduction}

Throughout, we will be working over the field $\C$ of complex numbers, and curves will be assumed to be projective.
\begin{defn}If $C$ is a smooth curve of genus $g\geq 2$, we say that a point $p$ of $C$ is a \emph{subcanonical point} if $K_C \isom \mathcal{O}_C\big((2g-2)p\big)$. \end{defn}
Equivalently, $p$ is subcanonical if there exists a holomorphic differential on $C$ which vanishes at $p$ to order $2g-2$ and nowhere else.  Let $\scl_g\subseteq \mathcal{M}_{g,1}$ be the locus of pointed curves $(C,p)$ where $p$ is a subcanonical point of $C$.

\subsection{Main results}  A point $p$ on a smooth curve $C$ is a \emph{Weierstrass point} if its associated set of \emph{Weierstrass gaps} 
\[ \{n\in \Z_{\geq 0}:~~h^0(C,\mathcal{O}_C(np))=h^0(C,\mathcal{O}_C((n-1)p))\}\]
is not equal to $\{1,2,\ldots, g\}$.  A point is a subcanonical if and only if $2g-1$ is a gap, so subcanonical points are Weierstrass points.  

It follows from results in \cite{KZ} that in genus $g\geq 4$, $\scl_g$ has three components $\scl_g^{\text{hyp}}$, $\scl_g^{\text{odd}}$ and $\scl_g^{\text{even}}$, corresponding to hyperelliptic curves and to the cases where the associated theta characteristic $\mathcal{O}_C\big((g-1)p\big)$ is odd or even.  In Section \ref{generalcase}, we determine the set of Weierstrass gaps associated to the general point of each of these components:

\begin{thm*}Let $g\geq 4$, then \begin{enumerate}\item a general point of $\sclhyp_g$ has Weierstrass gaps $\{1,3,5,\ldots,2g-5,2g-3,2g-1\}$,
\item a general point of $\sclodd_g$ has Weierstrass gaps $\{1,2,3,\ldots,g-2,g-1,2g-1\}$, and
\item a general point of $\scleven_g$ has Weierstrass gaps $\{1,2,3,\ldots,g-2,g,2g-1\}$.
\end{enumerate}
\end{thm*}

In Section \ref{lowgenus}, we show that for $g\leq 5$, these are the only possible Weierstrass gap sequences for subcanonical points, and in Section \ref{cyclic} we construct for $g\geq 6$ loci within $\scl_g^{\text{odd}}$ and $\scl_g^{\text{even}}$ consisting of subcanonical points with various other gap sequences as branch points of cyclic covers.

\subsection{Weierstrass points}  Let $C$ be a smooth curve of genus $g$. If $p$ is a point of $C$, the \emph{vanishing sequence} (of the complete canonical series $|K_C|$) at $p$ is the sequence $0=a^{K_C}_0(p)<a^{K_C}_1(p)<\ldots<a^{K_C}_{g-1}(p)\leq 2g-2$ of orders of vanishing of the holomorphic differentials at $p$, so that \[\left\{a^{K_C}_k(p)\right\}=\left\{v_p(\omega):~\omega\in H^0(C,K_C)\right\}.\]  Equivalently, we may consider the sequence $0=\alpha^{K_C}_0(p)\leq\alpha_1^{K_C}(p)\leq\ldots\leq\alpha_{g-1}^{K_C}(p)\leq~g-1$ defined by $\alpha_k^{K_C}(p)=a_k^{K_C}(p)-k$, which we call the \emph{ramification sequence} (of the canonical series) at $p$.  

Historically, Weierstrass points were defined using the \emph{Weierstrass gap sequence} at $p$.  This consists of, in increasing order, the positive integers $n$ for which there does not exist a meromorphic function on $C$ with pole divisor $np$.  One can check, using the Riemann-Roch theorem, that the set of Weierstrass gaps is just $\{a_i^{K_C}(p)+1\}$.  (One direction is easy: if there were a meromorphic function $f$ with pole divisor $(a_i+1)p$ and a holomorphic differential $\omega$ with $v_p(\omega)=a_i$, then the product $f\omega$ would be a meromorphic differential with only a single simple pole at $p$.)

One advantage of thinking in terms of the gap sequence is that the set of non-gaps $\N-\{a_i^{K_C}(p)+1\}$, the set of positive numbers $n$ for which there does exist a meromorphic function on $C$ with pole divisor $np$, forms a semigroup under addition (multiplying functions adds pole orders).  This provides a necessary (but not sufficient--for example, see \cite{EH87}) condition for a sequence of numbers  $0=a_0<a_1<\ldots<a_{g-1}\leq 2g-2$ to be the vanishing sequence for the canonical series at some point of some curve.

\par  At a general point of our curve $C$, there is no meromorphic function having a pole only at $p$ of order less than $g+1$, or equivalently the vanishing sequence at $p$ is the smallest possible vanishing sequence, $0,1,\ldots,g-1$, and correspondingly the ramification sequence is just $0,0,\ldots,0$.  We call a point $p$ with a larger ramification sequence $(\alpha_k^{K_C}(p))$ a \emph{Weierstrass point} of \emph{weight} $w(\alpha)=\sum_{k=0}^{g-1} \alpha_k $.  

A subcanonical point is then a Weierstrass point with $\alpha_{g-1}=g-1$.  For example, a hyperelliptic curve $C$ of genus $g$ has a degree $2$ map $\fun{\pi}{C}{\PR^1}$.  Each of the $2g+2$ ramification points of $\pi$ is in fact a Weierstrass point of $C$ with the maximal possible vanishing sequence $0,2,4,\ldots,2g-2$ and ramification sequence $0,1,2,\ldots,g-1$, and is thus an example of a subcanonical point.  Our main goal will be to try to determine what ramification sequences a subcanonical point can have and to study the stratification of $\scl_g$ by these sequences.

\par Every smooth curve has only finitely many Weierstrass points, the sum of whose weights is $g(g^2-1)$.  On a general curve, there are $g(g^2-1)$ distinct Weierstrass points, each of weight one (cf. \cite{ACGH} I.E).  For a given ramification sequence $\alpha=(\alpha_k)$, every component of the locally closed subset
\[ \mathcal{C}_\alpha = \{ (C,p) |~\alpha^{K_C}(p)=\alpha\}\subseteq \mathcal{M}_{g,1}\]
has codimension at most the weight of $\alpha$.  We say that $(C,p)$ is a \emph{dimensionally proper} Weierstrass point if $\mathcal{C}_{\alpha(p)}$ has codimension exactly $w(p)$ in a neighborhood of $(C,p)$.  It is known, for example (cf. \cite{EH87} for this and some slightly stronger results), that if $\alpha$ is any such sequence $0=\alpha_0\leq\alpha_1\leq\ldots\leq\alpha_{g-1}\leq~g-1$ with weight at most $g/2$, then $\mathcal{C}_\alpha$ contains dimensionally proper points. 
\subsection{Theta characteristics}

A \emph{theta characteristic} or \emph{spin structure} on a smooth curve $C$ of genus $g$ is a line bundle $L$ of degree $g-1$ on $C$ satisfying $L \otimes L \cong K_C$.  In \cite{C87}, Cornalba constructs a compactified moduli space of curves with theta characteristics.  The parity of the number $h^0(C,L)$ is constant in families of curves with theta characteristics, and a theta characteristic is defined to be \emph{odd} or \emph{even} as this number is.  

Any smooth curve has exactly $2^{g-1}(2^g-1)$ odd theta characteristics and $2^{g-1}(2^g+1)$ even theta characteristics.  Since a subcanonical point $p\in C$ is a point for which $\mathcal{O}_C((2g-2)p)\cong K_C$, the line bundle $L=\mathcal{O}_C ((g-1)p)$ is a theta characteristic, and we call the subcanonical point odd or even as this associated theta characteristic is. 
\par The parity of a subcanonical point $p\in C$ is in fact determined by its vanishing sequence: if $p$ is subcanonical, then $h^0\big(C,\mathcal{O}_C((g-1)p)\big)=h^0\big(C,K_C(-(g-1)p)\big)$ is the dimension of the space of holomorphic differentials vanishing to order at least $g-1$ at $p$, which is the number of $a_i^{K_C}(p)$ which are at least $g-1$.  This allows us to compute, for example that if $p$ is a ramification point of a hyperelliptic curve $C$, then $h^0(C,L)=\lfloor\frac{g+1}{2}\rfloor$, so $p$ is an odd subcanonical point if $g\equiv1,2\pmod{4}$ and even if $g\equiv0,3\pmod{4}$.

\subsection{Previous results}
In \cite{KZ}, Kontsevich and Zorich studied the moduli spaces $\mathcal{H}_g (k_1,k_2,\ldots,k_n)$ of pairs $(C,\omega)$ where $C$ is a compact Riemann surface of genus $g$ and $\omega$ is a holomorphic differential on $C$ with exactly $n$ zeroes, with multiplicities $k_1,\ldots,k_n$, where $k_1+\ldots + k_n=2g-2$.  They showed that these spaces are smooth (as orbifolds) of dimension $2g+n-1$ and classified their connected components.  The central case they considered was the case of $\mathcal{H}_g (2g-2)$, where the holomorphic differential $\omega$ has only a single zero of order $2g-2$ at a point $p$ which is then a subcanonical point of $C$.  They proved that for $g\geq 4$, the space $\mathcal{H}_g (2g-2)$ has three disjoint components, each of dimension $2g$, namely the locus where $C$ is hyperelliptic, and the loci where $C$ is non-hyperelliptic and $p$ is even and odd.  It follows immediately that  $\scl_g\subseteq \mathcal{M}_{g,1}$ has exactly three irreducible components, each of dimension $2g-1$, namely the hyperelliptic and the non-hyperelliptic even and odd subcanonical points. (The dimension is one lower than that in \cite{KZ} because of the freedom to multiply $\omega$ by a non-zero scalar without changing the point.)

\section*{Acknowledgements}
This paper is based on my Ph.D. dissertation at Harvard University.  I would like to thank my advisor, Joe Harris, for the very great deal of help he provided on this project.  This work was partially supported by an NSF Graduate Fellowship.

\bignote{Is it true that $\mathcal{H}(2g-2)$ is a $\C^{*}$-bundle over $\scl_g$ and what exactly does that tell us}

\section{General subcanonical points}\label{generalcase}
It follows from the results of \cite{KZ} that for $g\geq 4$, the locus $\scl_g\subseteq \mathcal{M}_{g,1}$ consisting of pairs $(C,p)$ where $p$ is a subcanonical point of $C$, has three disjoint irreducible components: \begin{enumerate}
\item the locus $\sclhyp_g$ of pairs $(C,p)$ where $C$ is a hyperelliptic curve of genus $g$ and $p$ is a ramification point of the hyperelliptic double cover, 
\item the locus $\sclodd_g$  of pairs $(C,p)$ where $C$ is a non-hyperelliptic curve of genus $g$ and $p$ is a subcanonical point such that $(g-1)p$ is an odd theta characteristic, and \bignote{nonhyperelliptic? non-hyperelliptic?}
\item the locus $\scleven_g$  of pairs $(C,p)$ where $C$ is a non-hyperelliptic curve of genus $g$ and $p$ is a subcanonical point such that $(g-1)p$ is an even theta characteristic.\end{enumerate}

In this section, we will prove the following theorem, which essentially states that a general subcanonical point $(C,p)$ in each component of $\scl_g$ has the smallest ramification sequence $\alpha^{K_C}(p)$ possible.

\begin{thm}\label{mainresult} Let $g\geq 4$, then \begin{enumerate}\item a general point of $\sclhyp_g$ has ramification sequence
$0,1,2,\ldots,g-3,g-2,g-1$,
\item a general point of $\sclodd_g$ has ramification sequence $0,0,0,\ldots,0,0,g-1$, and
\item a general point of $\scleven_g$ has ramification sequence $0,0,0,\ldots,0,1,g-1$.
\end{enumerate}
\end{thm}

In the hyperelliptic case, every ramification point of the hyperelliptic double cover has this ramification sequence and these $2g+2$ ramification points are all the Weierstrass points on the curve.  This is a standard result (see Chapter I of \cite{ACGH})  whose statement we include here for completeness.  It will also follow as a special case from our work on cyclic covers in Section \ref{cyclic}. 

As for the odd and even cases, the above ramification sequences are as small as possible: if $p\in C$ is any subcanonical point on a curve of genus $g$, then certainly \[\alpha^{K_C}(p)\geq (0,\ldots,0,0,g-1)\] and if $p$ is an even subcanonical point, then we must have $h^0(C,\mathcal{O}((g-1)p))\geq 2$; it follows that $a_{g-2}^{K_C}(p)\geq g-1$, so that \[\alpha^{K_C}(p)\geq (0,\ldots,0,1,g-1).\]

This means, by the upper semi-continuity of the ramification sequence, that in order to prove Theorem \ref{mainresult} we need only show that there exist points of $\sclodd_g$ and $\scleven_g$ having the desired ramification sequences.  

The basic approach of our construction will be that used in \cite{EH87} to construct Weierstrass points of low weight having prescribed ramification sequences: using limit linear series (cf. \cite{LLSBT}), we will begin by describing possible limit ``subcanonical points'' with the desired ramification sequences on certain reducible curves, and then show that they smooth to points with the same ramification sequences on nearby smooth curves.

We recall briefly the definition of a limit linear series in the case we will need.  Suppose that $X=C\cup E$ is the union of a smooth curve $C$ of genus $g-h$ and a smooth curve $E$ of genus $h$, meeting at a single node $q$.  Then a \emph{crude limit} $\grd{r}{d}$ on $X$ consists of a $\grd{r}{d}$ on $C$, i.e. a pair $L_C=(\mathcal{L}_C,V_C)$ where $\mathcal{L}_C$ is a line bundle of degree $d$ on $C$ and $V_C$ is a dimension $r+1$ subspace of $H^0(C,\mathcal{L}_C)$, together with a $\grd{r}{d}$ on $E$,  $L_E=(\mathcal{L}_E,V_E)$, satisfying the compatibility condition that for $0\leq i\leq r$,
\[a_i^{L_C}(q)+a_{r-i}^{L_E}(q) \geq d.\]

Given a one-parameter family of curves whose special fiber is $X$ and whose general fiber is a smooth curve of genus $g$, along with a family of $\grd{r}{d}$'s on the smooth fibers, there is a well-defined crude limit $\grd{r}{d}$ on $X$.  Strict inequalities in the compatibility condition arise when ramification points of the $\grd{r}{d}$'s approach the node of $X$; this complication can be avoided by blowing up the family at the node and base-changing.  A \emph{(refined) limit} $\grd{r}{d}$ on $X$ is a crude limit $\grd{r}{d}$ in which the above inequalities are all equalities.

The following lemma provides a partial description of the crude and refined limit canonical series on $X$ that have a ``subcanonical point'' on $E$ (i.e. a point $p$ with $a^{L_E}_{g-1}(p)=2g-2$) in the case where $E$ has genus $1$.

\begin{lem}\label{limitsc}Let $X=C\cup E$ be the union of a smooth curve $C$ of genus $g-1$ and a smooth curve $E$ of genus $1$, meeting at a single node $q$.  Assume $g\geq 3$.

\begin{enumerate} \item Let $L=\{(\mathcal{L}_C,V_C),(\mathcal{L}_E,V_E)\}$ be a crude limit $\grd{g-1}{2g-2}$ on $X$.  
\begin{enumerate} \item Suppose that there is a point $p\not=q$ of $E$ such that $\alpha_{g-1}^{L_E}(p)=g-1$.  Then $\mathcal{O}_E(p-q)\in \Pic^0(E)$ is $(2g-2)$-torsion, $q$ is a subcanonical point of $C$ (i.e. $\alpha_{g-2}^{K_C}(q)=g-2$), and $\alpha_0^{L_E}(p)=0$.
\item If moreover $\mathcal{O}_E(p-q)$ has order exactly $2g-2$ in $\Pic^0(E)$, then $\alpha_{i-1}^{K_C}(q)\geq  \alpha_{i}^{L_E}(p)$ for all $1\leq i \leq g-2$.
\item If instead  $\mathcal{O}_E(p-q)$ has order exactly $g-1$, then $\alpha_{i-1}^{K_C}(q)\geq  \alpha_{i}^{L_E}(p)$ for all but possibly one value of $i$ satisfying $1\leq i \leq g-2$, for which it must be the case that $g-1=a_i^{L_E}(p)=\alpha_i^{L_E}(p)+i$ and $\alpha_{i-1}^{K_C}(q)= \alpha_{i}^{L_E}(p)-1$.

\end{enumerate}
\item Suppose conversely that $q$ is a subcanonical point of $C$ with ramification sequence 
\[\alpha=(\alpha_0,\alpha_1,\ldots,\alpha_{g-3},g-2)\] and that $p$ is a point on $E$ such that $\mathcal{O}_E(p-q)$ has order $2g-2$ or $g-1$.  Then on $X$ there is a unique (refined) limit $\grd{g-1}{2g-2}$, $L=\{(\mathcal{L}_C,V_C),(\mathcal{L}_E,V_E)\}$ satisfying the ramification condition 
\[\alpha^{L_E}(p) \geq (0,\alpha_0,\alpha_1,\ldots,\alpha_{g-3},g-1)\]
at $p$.  Then $\alpha_{i+1}^{L_E}(p)=\alpha_{i}$ for $0\leq i \leq g-3$, except in the special case where $\mathcal{O}_E(p-q)$ has order $g-1$, $\alpha_i=g-3-i$, and either $i=g-3$ or $\alpha_{i+1}>\alpha_i$; in this case, which may only occur for at most one value of $i$, instead $\alpha_{i+1}^{L_E}(p)=\alpha_{i}+1$.
\end{enumerate}
\end{lem}

\begin{proof}
We begin by noting that if $L=\{(\mathcal{L}_C,V_C),(\mathcal{L}_E,V_E)\}$ is a crude limit $\grd{g-1}{2g-2}$ on $X$, then by ``Clifford's Theorem'' (Theorem 4.1 of \cite{LLSBT}) we know that 
$\mathcal{L}_C \isom K_C(2q)$ and that $\mathcal{L}_E\isom K_E(2(g-1)q) \isom \mathcal{O}_E((2g-2)q)$.  This completely determines the $C$-aspect of $L$: by Riemann-Roch we know that $h^0(C,K_C(2q))=g$, so that $V_C=H^0(C,K_C(2q))$, of which $H^0(C,K_C)$ is a codimension one subspace.  We see then that $a_0^{L_C}(q)=0$ and $a_i^{L_C}(q)=a_{i-1}^{K_C}(q)+2$ for $1\leq i \leq g-1$, or equivalently that $\alpha_0^{L_C}(q)=0$ and $\alpha_i^{L_C}(q)=\alpha_{i-1}^{K_C}(q)+1$  for $1\leq i \leq g-1$.  We thus need only consider the $E$-aspects of these linear series.

\par Now, to prove the first part of the lemma, we see that $\alpha_{g-1}^{L_E}(p)=g-1$ means $a_{g-1}^{L_E}(p)=2g-2$, so there exists some $\sigma\in V_E\subseteq H^0(E,\mathcal{L}_E)$ vanishing at $p$ to order $2g-2$, but $\mathcal{L}_E\isom\mathcal{O}_E((2g-2)q)$, so $(2g-2)p \sim (2g-2)q$ on $E$ and $\mathcal{O}_E(p-q)$ is $(2g-2)$-torsion as desired.  Moreover, as $\sigma$ does not vanish at $q$, we have $a_0^{L_E}(q)=0$, whence by the basic inequality for crude limit series \[a_0^{L_E}(q)+a_{g-1}^{L_C}(q)\geq 2g-2,\]
we conclude that $a_{g-2}^{K_C}(q) = a_{g-1}^{L_C}(q)-2 \geq 2g-2-2=2(g-1)-2$ and that $q$ is a subcanonical point of $C$.  \par Likewise, if $\alpha_0^{L_E}(p)>0$, then every $\sigma\in V_E$ would vanish at $p$ so that $a_{g-1}^{L_E}(q)< 2g-2$.  This, however, would imply that $a_0^{L_C}(q)>0$ by the basic inequality, which we know can not be the case.  Hence we must have $\alpha_0^{L_E}(p)=0$.

\par To show that if moreover $\mathcal{O}_E(p-q)$ has order exactly $2g-2$ then $\alpha_{i-1}^{K_C}(q)\geq  \alpha_{i}^{L_E}(p)$ for $1\leq i \leq g-2$, we use a simplified version of the proof of Proposition 5.2 of \cite{EH87}.  We would like to show that
\[a_{g-1-i}^{L_E}(q)+a_{i}^{L_E}(p)\leq 2g-3\mbox{~~~for~~} 1\leq i \leq g-2 \] 
since then by the basic inequality
 \[a_{g-1-i}^{L_E}(q)+a_{i}^{L_C}(q)\geq 2g-2\]
 we would know that $a_{i}^{L_C}(q)-1\geq a_{i}^{L_E}(p)$ and hence that 
 $\alpha_{i-1}^{K_C}(q)=\alpha_i^{L_C}(q)-1\geq \alpha_{i}^{L_E}(p)$ for $1\leq i \leq g-2$ as desired.
\par To prove our claim, suppose that $a_{g-1-i}^{L_E}(q)+a_{i}^{L_E}(p)\geq 2g-2$.  Let $W_1, W_2\subseteq V_E$ be the subspaces of $V_E$ defined by
\begin{align*}W_1&=\{\sigma\in V_E:~v_q(\sigma)\geq a_{g-1-i}^{L_E}(q)\}~\mbox{and}\\ 
                       W_2&=\{\sigma\in V_E:~v_p(\sigma)\geq a_{i}^{L_E}(p)\}.    \end{align*}
Then $\dim W_1=g-(g-1-i)=i+1$ and $\dim W_2=g-i$, and as these are subspaces of a $g$-dimensional vector space, we must have $\dim W_1\cap W_2\geq 1$.  Let $\sigma\in W_1\cap W_2$ be a nonzero element of the intersection.  We know that $\sigma$ has $2g-2$ zeroes, but since $a_{g-1-i}^{L_E}(q)+a_{i}^{L_E}(p)\geq 2g-2$ we have accounted for all of these, i.e. we must have that $a_{g-1-i}^{L_E}(q)+a_{i}^{L_E}(p)= 2g-2$ and
\[(2g-2)q\sim (\sigma)=a_{g-1-i}^{L_E}(q)q+a_{i}^{L_E}(p)p\]
so that $-(2g-2)q+a_{g-1-i}^{L_E}(q)q+a_{i}^{L_E}(p)p = a_{i}^{L_E}(p)(p-q) \sim 0$ and since $\mathcal{O}_E(p-q)$ has order $2g-2$, we must have $a_{i}^{L_E}(p)=0$ or $a_{i}^{L_E}(p)=2g-2$, so that $i=0$ or $i=g-1$.  
This shows that $a_{g-1-i}^{L_E}(q)+a_{i}^{L_E}(p)\leq 2g-3$  for $1\leq i \leq g-2$ as desired.
\par The case where $\mathcal{O}_E(p-q)$ has order exactly $g-1$ is entirely analogous, except that now 
\[(2g-2)p\sim (g-1)p+(g-1)q,\]
so it is possible that $a_{g-1-i}^{L_E}(q)+a_{i}^{L_E}(p)= 2g-2$ when $a_{g-1-i}^{L_E}(q)=a_{i}^{L_E}(p)=g-1$, in which case our use of the basic inequality as above only yields $\alpha_{i-1}^{K_C}(q)\geq \alpha_{i}^{L_E}(p)-1$.  Since the $a_i$ are increasing, this case may only occur for at most one value of $i$.  This completes the proof of the first part of the lemma.

\par
The second part of the lemma follows more or less directly from Proposition 5.2 of \cite{EH87}, which states in this case that given points $p$ and $q$ on an elliptic curve $E$ and given ramification sequences $\beta$ and $\gamma$ with 
\[\tag{*}g-1\geq \beta_{g-1-j}+\gamma_{j}\geq g-2\mbox{~~~for~~~} j=0,\ldots,g-1,\]
there exists at most one $\grd{g-1}{2g-2}$, $L=(\mathcal{O}_E((2g-2)q),V)$ on $E$ with $\alpha^L(q)=\beta$ and $\alpha^L(p)=\gamma$, and one exists if and only if 
\begin{equation*}\tag{**}\begin{split}\beta_{g-1-j}+\gamma_{j}=g-1&\implies b_{g-1-j}q+c_j p\sim (2g-2)q,  ~~~\mbox{and}\\ 
                       b_{g-1-j}q+(c_j+1)p \sim (2g-2)q &\implies \gamma_{j+1}=\gamma_j, \end{split}\end{equation*}
for each $j$, where $b_j=\beta_j+j$ and $c_j=\gamma_j+j$ are the associated vanishing sequences.
\par Now, the ramification sequence at $q$ is completely determined by the refined limit linear series condition and the known ramification sequence of the $C$-aspect at $q$: we must have
\[\beta=(0,g-2-\alpha_{g-3},g-2-\alpha_{g-2},\ldots,g-2-\alpha_{0},g-1)\]
so that $\beta_j=g-2-\alpha_{g-2-j}$ for $j=1,\ldots,g-2$.  On the other hand, we must consider all ramification sequences $\gamma$ at $p$ satisfying 
\[\gamma\geq\gamma'= (0,\alpha_0,\alpha_1,\ldots,\alpha_{g-3},g-1)\]
and show that for only one of them does there actually exist a $\grd{g-1}{2g-2}$, $L=(\mathcal{O}_E((2g-2)q),V)$ with  $\alpha^L(q)=\beta$ and $\alpha^L(p)=\gamma$; we will show that such a $\gamma$ must satisfy the hypothesis (*) of Proposition 5.2 of \cite{EH87}, so uniqueness will then be automatic.
\par In the case where $\mathcal{O}_E(p-q)$ has order $2g-2$, we must have $\gamma=\gamma'$, since if 
$\gamma_i>\gamma'_i=\alpha_{i-1}$ were bigger then the first part of the lemma would imply that $\alpha_{i-1}=\alpha^{K_C}_{i-1}(q)\geq \alpha_i^{L_E}(p)=\gamma_i>\alpha_{i-1}$.  In this case, condition (**) is clearly satisfied since $bq+cp\sim (2g-2)q$ is impossible except when $b$ and $c$ are $0$ and $2g-2$.
\par In the case where instead $\mathcal{O}_E(p-q)$ has order $g-1$, however, there might be a $\grd{g-1}{2g-2}$
with  $\alpha^L(q)=\beta$ and $\alpha^L(p)=\gamma\geq \gamma'$ where
$\gamma_i>\gamma'_i$ for some $i$.  However, the same argument using the first part of the lemma 
shows that this can happen for at most one index $i$ and that we must have $\gamma_i=\gamma'_i+1$ for that index.  It also shows that $i$ is completely determined by the $\alpha_j$: it must satisfy 
\[g-1=b_{g-1-i}= c'_i = c_i+1=\alpha_{i-1}+i+1.\]
\par When no such $i$ exists, our only possible ramification sequence at $p$ is $\gamma$, and Proposition 5.2 of \cite{EH87} tells us that there is a unique $\grd{g-1}{2g-2}$ on $E$ with the given ramification sequences at $p$ and $q$, since the second condition of (**) that 
\[b_{g-1-i}q+(c_i+1)p \sim (2g-2)q \implies \gamma_{i+1}=\gamma_i\]
is vacuous as there is no $i$ which makes $b_{g-1-i}=c_i+1=g-1$.
\par On the other hand, when there is such an $i$, we may have two different sequences of numbers, $\gamma$ and $\gamma'$, at $p$ which each satisfy the hypothesis (*) of the proposition together with the sequence $\beta$ at $q$, where again $\gamma$ and $\gamma'$ are equal except that $\gamma_i=\gamma'_i+1$.  Of course, in order for $\gamma'$ to actually be a ramification sequence as well, we must have that \[\gamma_i+1=\gamma'_i\leq \gamma'_{i+1}=\gamma_{i+1}\]  If this is not the case, then $\gamma_{i+1}=\gamma_i$, and hypothesis (**) is satisfied for $\beta$ and $\gamma$.
\par We are left with the case where there is such an $i$ but $\gamma_{i+1}\geq \gamma_i+1$.  In this case, $\beta$ and $\gamma$ do not satisfy hypothesis (**) since $\gamma_{i+1}\not=\gamma_i$, but $\beta$ and $\gamma'$ do satisfy hypothesis (**) since now $b_{g-1-j}q+c_j p\sim (2g-2)q$.  We have seen in each case then that there is a unique limit $\grd{g-1}{2g-2}$ on $X$ satisfying the given ramification condition at $p$, which completes the proof of the lemma.
\end{proof}

\subsection{Odd case}

We can now prove these exists a point of $\sclodd_g$ which has ramification sequence $0,\ldots,0,g-1$.  
As in \cite{EH87}, the proof is by induction on the genus.  We show in Section \ref{genus3} that the result is true in genus $3$, which provides the base case for our induction.

Suppose then that the result is known to be true in genus $g-1$, so that there exists a subcanonical point $q\in C$ on a curve of genus $g-1$ with ramification sequence $0,\ldots,0,g-2$.  As a Weierstrass point, $q$ is dimensionally proper: we know from \cite{KZ} that $\sclodd_{g-1}$ has dimension $2(g-1)-1 = 3(g-1)-2-(g-2)$.

As in \cite{EH87}, we consider the curve of compact type $X=C\cup E$ consisting of the curve $C$ meeting some elliptic curve $E$ in a node at $q$.  We pick a point $p\in E$ so that $\mathcal{O}_E(p-q)\in \Pic^0(E)$ has order $2g-2$.  Now, by Lemma \ref{limitsc}, there is a unique limit $\grd{g-1}{2g-2}$, $L=((\mathcal{L}_C,V_C),(\mathcal{L}_E,V_E))$ on $X$ with ramification sequence $0,\ldots,0,g-1$ at $p$.  Moreover, there is no such limit $\grd{g-1}{2g-2}$ satisfying the ramification condition at a general point $p'\in E$, since by the first part of Lemma \ref{limitsc}, $\mathcal{O}_E(p'-q)$ must be $(2g-2)$-torsion.  

Thus, as in Proposition 5.1 of \cite{EH87} we may apply Corollary 3.7 of \cite{LLSBT} to conclude that $L$ may be smoothed, maintaining the ramification condition $\alpha\geq (0,\ldots,0,g-1)$ near $p$.  We briefly sketch the Eisenbud-Harris argument:

Let $(\widetilde{X}\rightarrow \widetilde{B},\widetilde{B} \xrightarrow{\bar{p}} \widetilde{X})$ be a miniversal deformation space of the pointed curve $(X,p)$, with discriminant hypersurface $\Delta\subset \widetilde{B}$.  Then $\dim \widetilde{B} = 3g-2$.  By Theorem 3.3 of \cite{LLSBT}, there is a scheme \[G=G^{g-1}_{2g-2}\left(\widetilde{X}/\widetilde{B}; \left(\bar{p},\left(0,0,\ldots,0,g-1\right)\right)\right)\longrightarrow \tilde{B}\]
whose fiber over each point $z\in \tilde{B}$ parametrizes the (refined) limit $\grd{g-1}{2g-2}$'s $L_z$ on $\widetilde{X}_z$ satisfying the ramification condition
\[\alpha^{L_z}(\bar{p}(z))\geq (0,0,\ldots,0, g-1)\]
at the marked point $\bar{p}(z)$ of $\widetilde{X}_z$.  Moreover, every component of $G$ has dimension at least the expected dimension:
\begin{align*}
\dim G &\geq \dim \widetilde{B} + \rho(g,g-1,2g-2; (0,0,\ldots,0,g-1)) \\
& = \dim \widetilde{B} + (g-1+1)(2g-2-(g-1))-(g-1)g-(g-1)\\
& = (3g-2)  - (g-1)\\
&= 2g-1.
\end{align*}
Suppose that the component of $G$ containing our given limit $\grd{g-1}{2g-2}$ $L$ on $(X,p)$ were to lie entirely over $\Delta$ (i.e. suppose that $L$ doesn't smooth).  Then by the first part of Lemma \ref{limitsc}, $G$ must in fact lie over the locus $D\subset \Delta$ parameterizing pointed nodal curves $C_z\cup_{\bar{q}(z)} E_z$ in which $\bar{q}(z)$ is a subcanonical point of $C_z$ and $\mathcal{O}_{E_z}(\bar{p}(z)-\bar{q}(z))$ is $(2g-2)$-torsion.  Since, by the induction hypothesis, we know that the locus of genus $g-1$ curves $(C_z,\bar{q}(z))$ with a marked subcanonical point has dimension $(3(g-1)-2)-((g-1)-1)=2g-3$, and since the locus of elliptic curves with a marked point of order $(2g-2)$ has dimension $1$, we find that $\dim D = 2g-2$.  But we know by the second part of Lemma \ref{limitsc} that the curves over $D$ each have a unique $\grd{g-1}{2g-2}$ satisfying the ramification condition, i.e. that the fibers of $G$ over $D$ each consist of a single point.  Thus the dimension of this component of $G$ must be $2g-2$, contradicting the known bound $\dim G\geq 2g-1$.  We conclude that $G$ does not lie entirely over $\Delta$.

This means there are nearby smooth pointed curves with ramification sequence (of their canonical series, the only $\grd{g-1}{2g-2}$ on a smooth curve) at least $0,\ldots,0,g-1$.  However, $L$ itself has exactly this ramification sequence at $p$, so by upper semi-continuity of ramification sequences, the general smooth curve in the smoothing family must also have ramification sequence exactly $0,\ldots,0,g-1$, as desired.  This completes our proof by induction that a general point of $\sclodd_g$ has ramification sequence $0,\ldots,0,g-1$.

\subsection{Even case}

We turn now to the case of $\scleven_g$.  While we know from \cite{KZ} that $\scleven_g$ has the same dimension, $2g-1$, as $\sclodd_g$, the smallest possible ramification sequence $0,\ldots,0,1,g-1$ has a weight which is one greater.  This means that a general point of $\scleven_g$ cannot possibly be dimensionally proper, which is a problem: the inductive framework of \cite{EH87} only applies to dimensionally proper Weierstrass points.

\par In order to deal with this issue, we note that while these are not dimensionally proper Weirstrass points, if we think of $\scleven_g$ as being the subcanonical points whose associated theta-characteristic is even, then it does have the ``expected dimension'' in the sense that while the ramification condition $\alpha^{K_C}_{g-1}(p)\geq g-1$ ought give a codimension of $g-1$, we might expect that the additional condition that $\mathcal{O}_C\big((g-1)p\big)$ be an even theta-characteristic would not increase the codimension further since the parity of a theta-characteristic is constant in families.

For the even case, we consider the same nodal curve $X=C\cup E$ as in the odd case, where an elliptic curve $E$ meets a curve $C$ of genus $g-1$ at an \emph{odd} subcanonical point $q$ of $C$ with ramification sequence $0,\ldots,0,g-2$.  Now, we pick a distinguished point $p$ of $E$ such that $\mathcal{O}_E(p-q)$ has order exactly $g-1$, instead of $2g-2$ as in the even case.  We will see shortly that this corresponds to picking an odd theta characteristic on $E$ rather than an even one.

As in the odd case, the point of attachment is a dimensionally proper Weierstrass point and Lemma \ref{limitsc} again shows that there is a unique limit $\grd{g-1}{2g-2}$, $L=((\mathcal{L}_C,V_C),(\mathcal{L}_E,V_E))$ on $X$ satisfying the same ramification condition \[\alpha^{L_E}(p)\geq (0,\ldots,0,g-1)\] at $p$ and that there is no such limit $\grd{g-1}{2g-2}$ satisfying the ramification condition at a general point of $E$.  Thus, the proofs of Proposition 5.1 of \cite{EH87} and Corollary 3.7 of \cite{LLSBT} again may be applied to show that $L$ can be smoothed, preserving the ramification condition $\alpha \geq (0,\ldots,0,g-1)$ near $p$.

\par More precisely, there is a family of stable curves $\fun{\pi}{\mathcal{X}}{\Delta}$, over a smooth, one-dimensional base, with smooth fibers $X_t=\pi^{-1}(t)$ away from the special fiber $X=X_0=C\cup E$, together with a section $\fun{\bar{p}}{\Delta}{\mathcal{X}}$ such that $\bar{p}(0)=p$ and 
\[\alpha^{K_{X_t}} (\bar{p}(t)) \geq 0,\ldots,0,0,g-1\]
for $t\not=0$, as again the complete series $K_{X_t}$ is the only $\grd{g-1}{2g-2}$ on the smooth curve $X_t$.

We thus know that the points $\bar{p}(t)$ are subcanonical points of the smooth curves $X_t$ for $t\not=0$, and we know by Lemma \ref{limitsc} that $\alpha^{L_E}(p)=0,\ldots,0,1,g-1$.  By upper semi-continuity, this leaves only two possibilities for the ramification sequence at $\bar{p}(t)$ for a general $t$, namely $0,\ldots,0,0,g-1$ and $0,\ldots,0,1,g-1$.  In order to show that it is in fact the latter, we will need to apply a result about the limits of theta characteristics on smooth curves approaching a curve of compact type.

In \cite{C87}, Cornalba constructs a compactified moduli space of curves with theta-characteristics
by describing objects associated to a stable curve which correspond to limits of theta-characteristics on nearby smooth curves.  In this compactificaton, the odd and even loci remain disjoint irreducible components (cf. section 6 of \cite{C87}).
\par In the case of curves of compact type, the answer is especially simple: a ``theta-characteristic'' on a curve of compact type should consist of a choice of a theta-characteristic on each of its components of positive genus, and its parity should be the sum of the parities on those components.  The special case that we require is the following: 
\begin{lem}\label{thetaparity}Let $\fun{\pi}{\mathcal{X}}{\Delta}$ be a family of stable curves of genus $g$ over a smooth, one-dimensional base, such that for some $0\in \Delta$, the family is smooth away from $0$ and the special fiber $\pi^{-1}(0)=C\cup C'$ consists of a curve $C$ of genus $i$ and a curve $C'$ of genus $g-i$ meeting at a single node.
\par Let $\mathcal{L}_0\rightarrow(\mathcal{X}-\pi^{-1}(0))$ be a family of theta-characteristics on the smooth fibers of $\mathcal{X}$.  Let $\mathcal{L}$ be the unique extension of $\mathcal{L}_0$ to all of $\mathcal{X}$ which has degree $i-1$ on $C$ and let $\mathcal{L}'$ be the unique extension of $\mathcal{L}_0$ to all of $\mathcal{X}$ which has degree $g-i-1$ on $C'$.  
\par Then $\mathcal{L}|_C$ is a theta-characteristic on $C$, $\mathcal{L}'|_{C'}$ is a theta-characteristic on $C'$, and the parity of the theta characteristic $\mathcal{L}_0|_{\pi^{-1}(\lambda)}$ for $\lambda\not=0$ is equal to the sum of the parities of $\mathcal{L}|_C$ and $\mathcal{L}'|_{C'}$.  \end{lem}

In the case of our family $\mathcal{X}\to \Delta$, on each smooth curve $X_t$, the associated theta-characteristic is $\mathcal{O}_{X_t}((g-1)\bar{p}(t))$.  Thus one extension of this family of theta-characteristics to a line bundle on all of $\mathcal{X}$ is simply $\mathcal{O}_\mathcal{X} ((g-1)\bar{p})$.  The theta characteristics on $C$ and $E$ associated to this family are thus simply the twists of this line bundle of degrees $g-2$ on $C$ and degree $0$ on $E$, respectively.  These are $\mathcal{O}_C((g-2)q)$ and $\mathcal{O}_E((g-1)p-(g-1)q)$.  Now, we know that $q$ is an odd subcanonical point of $C$ (since it has ramification sequence $0,\ldots,0,g-2$) so $\mathcal{O}_C((g-2)q)$ is an odd theta-characteristic.  On the other hand $\mathcal{O}_E((g-1)p-(g-1)q)\isom \mathcal{O}_E$ is effective, since $p$ was chosen so that $\mathcal{O}_E(p-q)$ has order $g-1$ in $\Pic^0(E)$, so the theta-characteristic on the elliptic curve $E$ is odd.  

By Lemma \ref{thetaparity}, $\mathcal{O}_{X_t}((g-1)\bar{p}(t))$ is an even theta-characteristic for $t\not=0$.  This implies that the ramification sequence of the subcanonical point $\bar{p}(t)$ for general $t$ is $0,\ldots,0,1,g-1$ rather than $0,\ldots,0,0,g-1$, completing the proof that the general point of $\scleven_g$ has ramification sequence $0,\ldots,0,1,g-1$.

\vspace{14pt}

\begin{Rem} A slightly closer examination of the proof of Theorem \ref{mainresult}, in particular of the use of Corollary 3.7 of \cite{LLSBT}, also provides a new proof, without using the methods of \cite{KZ}, that 
$\sclodd_g$ and $\scleven_g$ are non-empty for $g\geq 4$ and that each has some component of the correct dimension $2g-1$.  It seems likely that this proof, unlike that in \cite{KZ}, might extend to the case of characteristic $p$, at least when $p\gg g$, using the theory of limit linear series in characteristic $p$ developed in \cite{charplls}.

There does not, however, seem to be any easy way to show the irreducibility of $\sclodd_g$ and $\scleven_g$ using these techniques.
\end{Rem}

\begin{Rem}
 While the proof of Theorem \ref{mainresult} would generalize using any other dimensionally proper subcanonical point as a base case for the induction, in fact no such points can exist.  The corollary to Theorem 2 of \cite{EH87} shows that a dimensionally proper point $p\in C$ must be \emph{primitive}, which means that all smaller sequences $\alpha \leq \alpha^{K_C}(p)$ must satisfy the semigroup condition.  One can check directly that the sequence 
 \[0,0,\ldots,0,1,g-2\]
 fails to satisfy the semigroup condition and that thus $0,0,\ldots,0,g-1$ is the only possible ramification sequence for a dimensionally proper subcanonical point. 
 
 \end{Rem}

 \section{Cyclic covers}\label{cyclic}
 Perhaps the simplest construction of subcanonical points is as ramification points of the hyperelliptic double cover on a hyperelliptic curve.  In this section, we show how non-hyperelliptic subcanonical points can be constructed as ramification points of certain cyclic covers, and describe in some cases how to compute the vanishing sequences of those subcanonical points.
 
 Let $d>1$ be a fixed natural number, and $B$ be a given smooth curve of genus $h$.  Let $D$ be an effective divisor on $B$ consisting of distinct points, with $\deg D$ divisible by $d$.  Let $L$ be a line bundle on $B$ satisfying $L^{\otimes d}\isom \mathcal{O}_B(D)$.  Then we may construct a $d$-sheeted cyclic cover of $B$, totally ramified over each point of $D$, as follows:
the isomorphism $L^{\otimes d}\isom \mathcal{O}_B(D)$ determines a $d$th power map from the total space of the line bundle $L$ to the total space of the line bundle $\mathcal{O}_B(D)$.  The total space of the line bundle $\mathcal{O}_B(D)$ has a distinguished section, the constant section $1$, which as a section of $\mathcal{O}_B(D)$ vanishes on $D$.  The preimage of this section in the total space on $L$ is the desired $d$-sheeted cyclic cover $\fun{\pi}{C}{B}$, totally ramified over $D$: it is easily checked that it has a local analytic equation of the form $z\mapsto z^n$ at the ramification points, and the cyclic automorphism group is giving simply by the group of $d$th roots of unity acting on the total space of $L$ by multiplication on each fiber.  Let $g$ be the genus of $C$.

We now give a more algebraic description of this construction (see \cite{hart} II.5, IV.3).  Consider the locally free sheaf of rank $d$
\[\mathcal{F}=\mathcal{O}_B\oplus L^{-1} \oplus L^{\otimes(-2)}\oplus\cdots\oplus L^{\otimes(-(d-1))}\] 
on $B$.  The maps 
$L^{\otimes(-i)}\otimes L^{\otimes(-j)} \isom L^{\otimes(-i-j)}$ and $L^{\otimes(-i)}\otimes L^{\otimes(-j)} \isom L^{\otimes(-i-j)} \isom L^{\otimes(-i-j+d)}(-D)\hookrightarrow L^{\otimes(-i-j+d)}$ give $\mathcal{F}$ the structure of a sheaf of $\Z / d \Z$-graded algebras over $\mathcal{O}_B$.  The desired cyclic cover is then the global spec
\[\fun{\pi}{C=\SPEC \mathcal{F}}{B},\]
and here the action of a $d$th root of unity $\zeta$ is induced by multiplication by $\zeta^i$ on the $i$th graded piece.
This tells us in particular that 
\[\pi_{*} \mathcal{O}_C \isom \mathcal{O}_B\oplus L^{-1} \oplus L^{\otimes(-2)}\oplus\cdots\oplus L^{\otimes(-(d-1))},  \]
and moreover, looking in local analytic coordinates near a ramification point $p\in C$ with $\pi(p)=q$ so that $\pi$ is given near $p$ by $z\mapsto z^n$, we see that the $L^{\otimes(-i)}$ component in this direct sum decomposition corresponds to functions which are multiplied by $\zeta^i$ when $z$ is replaced by $\zeta z$.  Equivalently, a holomorphic section of $L^{-i}$ near $q$ corresponds in this decomposition to a holomorphic function on a $z$-disc near $p$ in whose Taylor expansion $z^j$ may only appear when $j\equiv i\pmod{d}$.  This shows in particular that a section of $L^{\otimes(-i)}$ vanishing to order $m$ at $q$ corresponds to a holomorphic function on $C$ near $p$ which vanishes to order $dm+i$, and likewise for meromorphic sections and pole orders in the case where $m$ is negative.

Thus if $p$ is a ramification point of $\pi$ and $\pi(p)=q$, then
\[\pi_{*} \mathcal{O}_C(kp) \isom \mathcal{O}_B\left(\floor{\tfrac{k}{d}}q\right)  \oplus L^{-1}\left(\floor{\tfrac{k+1}{d}}q\right) \oplus L^{\otimes(-2)}\left(\floor{\tfrac{k+2}{d}}q\right)\oplus\cdots\oplus L^{\otimes(-(d-1))}\left(\floor{\tfrac{k+d-1}{d}}q\right).\]
We would like to determine when $p$ is subcanonical and calculate its vanishing sequence when it is subcanonical, or in other words we would like to compute the quantity 
\[h^0\big(C,K_C(-np)\big)=h^0\big(B,\pi_{*}K_C(-np)\big),\]
especially in the case where $n=2g-2$, where it will be one if $p$ is subcanonical and zero otherwise.  To do this, we note that by Riemann-Hurwitz, $\pi^{*}K_B\isom K_C(-R)$, where $R$ is the ramification divisor of $\pi$, consisting of each ramification point with multiplicity $d-1$, so that $\pi(R)=(d-1)D$.  We may then compute
\begin{align*}
\pi_{*}K_C(-np) \isom {} & \pi_* (\mathcal{O}_C(R-np)\otimes\pi^* K_B)  \\
                         \isom {} &\pi_* (\mathcal{O}_C(R-np))\otimes K_B \\
                         \isom {} & K_B\left(\floor{\tfrac{d-1-n}{d}}q\right) \oplus \left(K_B\otimes L^{-1}\left(D+\floor{\tfrac{-n}{d}}q\right)\right)\oplus \left(K_B\otimes L^{\otimes(-2)}\left(D+\floor{\tfrac{1-n}{d}}q\right)\right)\\ & \oplus\cdots\oplus \left(K_B\otimes L^{\otimes(-(d-1))}\left(D+\floor{\tfrac{d-2-n}{d}}q\right)\right)\\
                         \isom {} & K_B\left(\floor{\tfrac{d-1-n}{d}}q\right) \oplus \left(K_B\otimes L^{\otimes(d-1)}\left(\floor{\tfrac{-n}{d}}q\right)\right)\oplus \left(K_B\otimes L^{\otimes(d-2)}\left(\floor{\tfrac{1-n}{d}}q\right)\right)\\ & \oplus\cdots\oplus\left(K_B\otimes L^{\otimes 2}\left(\floor{\tfrac{d-3-n}{d}}q\right)\right)\oplus \left(K_B\otimes L\left(\floor{\tfrac{d-2-n}{d}}q\right)\right)     \\
                         \isom {} & \bigoplus_{i=0}^{d-1} L^{\otimes i}\otimes K_B \left(\floor{\tfrac{d-1-i-n}{d}}q\right),              
\end{align*}
since $L^{\otimes d}\isom \mathcal{O}_B(D)$.  By Riemann-Hurwitz, \[2g-2=d(2h-2)+(d-1)\deg D=d((2h-2)+(d-1)\deg L)\] and in particular $2g-2$ is a multiple of $d$.  When we set $n=2g-2$, only one of the line bundles in the above direct sum has non-negative degree, namely $L^{\otimes (d-1)}\otimes K_B\big(-\tfrac{2g-2}{d}q\big)$, which has degree $(d-1)\tfrac{2g-2-d(2h-2)}{d(d-1)}+(2h-2)-\tfrac{2g-2}{d}=0$.  We see then that 
\begin{align*}
\mbox{$p$ is a subcanonical point of $C$} & \Longleftrightarrow L^{\otimes(d-1)}\otimes K_B\isom \mathcal{O}_B\left(\tfrac{2g-2}{d}q\right) \\
& \Longleftrightarrow K_B(D)\isom L\left(\tfrac{2g-2}{d}q\right)\\
& \Longleftrightarrow L\isom K_B\left(D-\tfrac{2g-2}{d}q\right).
\end{align*}

Thus, given $q$ and $B$, we would like to determine whether there is some effective divisor $D$ on $B$, consisting of distinct points, one of which is $q$, so that 
\[\left(K_B\left(D-\tfrac{2g-2}{d}q\right)\right)^{\otimes d}\isom\mathcal{O}_B(D),\]
or equivalently, so that
\[(d-1)D\sim (2g-2)q-d K_B.\]
Now if we let $D'$ be some divisor (not necessarily effective) such that $(d-1)D'\sim (2g-2)q-d K_B$, we would like to find an effective divisor in $|D'-q|$ which consists of $\deg D' -1$ distinct points other than $q$.  By Bertini's theorem, this is possible as long as $|D'-q|$ is base-point-free, i.e. if 
\[0<h^0(B,\mathcal{O}_B(D'-q-r)) < h^0(B,\mathcal{O}_B(D'-q))\] 
for every $r\in B$.  By Riemann-Roch, this is guaranteed to be the case if $\deg D'-2\geq 2h-1$, i.e. if the number of distinct branch points of the cover we are constructing is at least $2h+1$.

We would now like to find the vanishing sequence of the subcanonical point $p\in C$ that we have constructed, for which we must compute $h^0\big(C,K_C(-np)\big)=h^0\big(B,\pi_{*}K_C(-np)\big)$ for $0\leq n\leq 2g-2$.  However, in our direct sum decomposition for $\pi_{*}K_C(-np)$, without additional hypotheses we only have control over the terms 
\[K_B\left(\floor{\tfrac{d-1-n}{d}}q\right)\mbox{~~and~~}K_B\otimes L^{\otimes(d-1)}\left(\floor{\tfrac{-n}{d}}q\right)\isom \mathcal{O}_B\left(\left(\tfrac{2g-2}{d}+\floor{\tfrac{-n}{d}}\right)q\right)\] which together determine only the portion of the vanishing sequence where $n\equiv 0,d-1 \pmod{d}$; for $n$ is in the vanishing sequence for $K_C$ at $p$ if and only if $h^0\big(C,K_C(-np)\big)\not=h^0\big(C,K_C(-(n+1)p)\big)$, and 
$\floor{\tfrac{i-n}{d}}\not=\floor{\tfrac{i-(n+1)}{d}}$ if and only if $n\equiv i\pmod{d}$, so whether a number $n\equiv i\pmod{d}$ appears in the vanishing sequence depends entirely on the $L^{\otimes(d-1-i)}\otimes K_B\left(\floor{\tfrac{i-n}{d}}\right)$ term.

We see then that $n=dm+(d-1)$ is in the vanishing sequence for $K_C$ at $p$ if and only if $m$ is in the vanishing sequence for $K_B$ at $q$.  

Likewise, if $n=dm$, then $n\in\{a^{K_C}_j(p)\}$ if and only if $h^0\left(B,\mathcal{O}_B\left(\left(\tfrac{2g-2}{d}-m\right)q\right)\right)\not=h^0\left(B,\mathcal{O}_B\left(\left(\tfrac{2g-2}{d}-(m+1)\right)q\right)\right)$, or in other words, if and only if $\tfrac{2g-2}{d}-m$ is a Weierstrass non-gap of $p$, or if and only if $0\leq m\leq\tfrac{2g-2}{d}$ and $\tfrac{2g-2}{d}-m-1\not\in\{a^{K_B}_j(q)\}$.

\subsection{Double covers}The other terms may depend on tensor powers of $K_B$ in ways that are not completely determined by the vanishing sequence of $K_B$ itself at $q$, so if we are to compute the entire vanishing sequence we will need to make an additional assumption.  One additional assumption we may make is that $d=2$, in which case the two terms we can control are the only terms.  This gives us the following:

\begin{thm}\label{doublecover}Let $B$ be a curve of genus $h$ and $q$ be a point of $B$.  Let $g\geq 3h$.  Then there exists a double cover $\fun{\pi}{C}{B}$, where $C$ has genus $g$, such that $q$ is a branch point of $\pi$, and the point $p\in C$ with $\pi(p)=q$ is a subcanonical point with vanishing sequence as follows:
\begin{align*}
2m+1\in \{a_j^{K_C}(p)\} &~~\Longleftrightarrow& m &\in \{a_j^{K_B}(q)\} \\
2m \in \{a_j^{K_C}(p)\} &~~\Longleftrightarrow& g-2-m &\not\in \{a_j^{K_B}(q)\} \text{ and } 0\leq m\leq g-1.
\end{align*}
Conversely, every subcanonical point on a curve of genus $g$ which is a branch point of a double cover to a curve of genus $h$ has vanishing sequence as above.
\end{thm}
The converse here follows from the fact that in the special case of $d=2$, every double cover can be constructed in the way we have described. When $d>2$ it is no longer the case that every $d$-sheeted cyclic cover arises in this way.  Applying Theorem \ref{doublecover} in the case $h=0$ recovers the standard computation of the vanishing sequence of a hyperelliptic Weierstrass point:
\begin{cor}\label{hyperelliptic}A branch point of the hyperelliptic double cover from a hyperelliptic curve $C$ of genus $g$ is a subcanonical point with vanishing sequence 
\[ 0, 2, 4,\ldots, 2g-4, 2g-2\]
and ramification sequence
\[ 0,1,2,\ldots,g-2,g-1.\]
\end{cor}
Applying Theorem \ref{doublecover} in the case $h=1$ recovers the computation of the vanishing sequence of a subcanonical point which is a ramification point on a bielliptic curve.  This is proven in \cite{bielliptic} in the course of classifying Weierstrass points on bielliptic curves.  Their results imply also that a subcanonical point of a bielliptic curve $\fun{\pi}{C}{E}$ which is not a ramification point of $\pi$ must have ramification sequence $0,\ldots,0,0,g-1$ or $0,\ldots,0,1,g-1$.
\begin{cor}\label{bielliptic}For $g\geq 3$, there exists a subcanonical point on a bielliptic curve $C$ of genus $g$ which is a branch point of the bielliptic double cover $C\to E$.  Any subcanonical point on a bielliptic curve which is the branch point of the bielliptic cover has vanishing sequence
\[ 0, 1, 2, 4, 6, \ldots, 2g-8, 2g-6, 2g-2\]
and ramification sequence
\[ 0, 0, 0,1,2 \ldots,g-5,g-4,g-1.\]
\end{cor}

In the case $h=2$, there are two cases for $q$: the point can be a general point or a Weierstrass point with ramification sequence $0,2$.  Applying Theorem \ref{doublecover} to those two cases we get:
\begin{cor}\label{doublecoversofgenus2} For $g\geq 6$, there exist subcanonical points on curves of genus $g$ with the ramification sequence
\[0, 0, 0, 0, 0, 1, 2 \ldots,g-7,g-6,g-1.\]
For $g=6$, there exist subcanonical points with ramification sequence $0,0,0,2,2,5$, for $g=7$, there exist subcanonical points with ramification sequence $0,0,0,1,1,3,6$, and for $g\geq 8$, there exist subcanonical points with ramification sequence
\[0, 0, 0, 1, 1, 1, 2, \ldots,g-8,g-7, g-4, g-1.\]
\end{cor}
In the case $h=3$, there are four different possible vanishing sequences for a point on a genus $3$ curve, namely $0,1,2$; $0,1,3$; $0,1,4$; and $0,2,4$.  Applying Theorem \ref{doublecover} to those four cases we get:
\begin{cor} For $g\geq 14$, there exist subcanonical points on curves of genus $g$ that have each of the following ramification sequences:
\begin{alignat*}{8}
&0,0,0,0,0,0,0,1,2,3,4,&\ldots,~&&&g-8,~&g-1&,\\
&0,0,0,0,0,1,1,1,2,3,4,&\ldots,~&&g-9,~&g-6,~&g-1&,\\
&0,0,0,0,0,1,2,2,2,3,4,&\ldots,~&g-10,~&g-7,~&g-6,~&g-1&\text{, and}\\
&0,0,0,1,1,1,2,2,2,3,4,~&\ldots,~&g-10,~&g-7,~&g-4,~&g-1&.
\end{alignat*}
\end{cor}
While this the proof of the corollary works for all $g\geq 9$, writing down the sequences it yields for $9\leq g \leq 13$ requires some care, as in the $g=6, 7$ cases of Corollary \ref{doublecoversofgenus2}.

   We could also let $h$ vary; for example, setting $h=\left\lfloor\tfrac{g}{3}\right\rfloor$, the biggest $h$ can be in Theorem \ref{doublecover} and taking $q$ to be a general point on a curve of genus $h$, we get:
\begin{cor} For $g\geq 6$, let $g\equiv i \pmod{3}$, with $i\in\{0,1,2\}$.  Then
there exist subcanonical points on curves of genus $g$ with ramification sequence
\[0,\ldots,0,1,2,3,\ldots,\left\lfloor\tfrac{g}{3}\right\rfloor-2+i,g-1.\]
\end{cor}

\subsection{Cyclic covers of higher degree}

Alternatively, instead of requiring $d=2$, we can allow $d>2$, but control the $L^{\otimes(d-1-i)}\otimes K_B\left(\floor{\tfrac{i-n}{d}}\right)$ terms by assuming that $(2h-2)q\sim K_B$; this is automatic in the cases where $h=0$ or $h=1$ and equivalent to saying that $q$ is itself a subcanonical point of $B$ when $h\geq 2$.  This simplifies the condition $(d-1)D\sim (2g-2)q-dK_B$ to $(d-1)D\sim (2g-2-d(2h-2))q$, and allows us to pick $D$ with 
\[D\sim \tfrac{2g-2-d(2h-2)}{d-1}q\mbox{~~and~~}L\isom \mathcal{O}_B\left(\tfrac{2g-2-d(2h-2)}{d(d-1)}q  \right).\]
 For the rest of the section, we will assume that $D$, and thus also $L$, is as above.  Letting $\ell=\deg L=\tfrac{2g-2-d(2h-2)}{d(d-1)}$, this yields
\[L^{\otimes(d-1-i)}\otimes K_B\left(\floor{\tfrac{i-n}{d}}\right)\isom K_B\left(\left(\floor{\tfrac{i-n}{d}}+(d-1-i)\ell\right)q\right),\]
and we see that $n=dm+i\in \{a_j^{K_C}(p)\}$ if and only if either $(d-1-i)\ell-m\leq 0$ with $m-(d-1-i)\ell \in \{a_j^{K_B}(q)\}$, or $d-1-i-m\geq 2$, since meromorphic $1$-forms exist with pole locus $kq$ for $k\geq 2$ but not for $k=1$.  Solving for $g$ in terms of $\ell$, this gives us:
\begin{thm}Let $B$ be a curve of genus $h$ and let $q$ be a point of $B$ such that $K_B\isom\mathcal{O}_B\big((2h-2)q\big)$. Let $d\geq 2$ be an integer and $\ell$ be a positive integer satisfying $\ell\geq \tfrac{2h+1}{d}$.  Then if $g=\tfrac{d(d-1)}{2}\ell+d(h-1)+1$ is at least $2$, there is a curve $C$ of genus $g$, which is a cyclic $d$-sheeted cover $\fun{\pi}{C}{B}$ of $B$ such that there is a totally ramified $p\in C$ with $\pi(p)=q$, where $p$ is a subcanonical point of $C$ whose vanishing sequence is determined by the following, for $i=0,1,\ldots,d-1$:
\[dm+i\in \{a_j^{K_C}(p)\} \Longleftrightarrow \begin{aligned}\text{ either~~}&0\leq m\leq(d-1-i)\ell-2\\\text{ or~~}&m-(d-1-i)\ell \in \{a_j^{K_B}(q)\}.\end{aligned}\]
\end{thm}

Taking $h=0$ or $h=1$, with $g=3\ell-2$ or $g=3\ell+1$, respectively, we get the following:
\begin{cor}
Let $g=3k+1\geq 7$.  Then there exist subcanonical points on curves of genus $g$ with each of the following ramification sequences:
\[0,0,1,1,2,2,\ldots, k-1, k-1, k, k+2, k+4, \ldots, g-5, g-3, g-1,\text{ and}\]
\[0,0,0,0,0,1,1,2,2,\ldots, k-3, k-3, k-2, k, k, k+1, k+3, \ldots, g-8, g-6, g-1.\]
\end{cor}

 \section{Low genus examples}\label{lowgenus}
 The following table summarizes the possible ramification sequences of subcanonical points on curves of genus $g\leq 6$, which we will be describing for the rest of the section.
 
 \begin{center}
 \begin{tabular}{|c|r|r|r|r|r|}
 \hline
 genus & $a^{K_C}(p)$ & $\alpha^{K_C}(p)$ & parity & weight & $\codim{M_{g,1}} \mathcal{C}_\alpha$\\ \hline
 $2$ & $0, 2$ & $0, 1$ & odd & $1$ & $1$\\ \hline
 
  $3$ & $0, 1, 4$ & $0, 0, 2$ & odd & $2$ & $2$\\ 
   & $0, 2, 4$ & $0, 1, 2$ & even & $3$ & $2$\\ \hline
 
  $4$ & $0, 1, 2, 6$ & $0, 0, 0, 3$ & odd & $3$ & $3$\\ 
   & $0, 1, 3, 6$ & $0, 0, 1, 3$ & even & $4$ & $3$\\ 
   & $0, 2, 4, 6$ & $0, 1, 2, 3$ & even & $6$ & $3$\\ \hline
   
$5$ & $0, 1, 2, 3, 8$ & $0, 0, 0, 0, 4$ & odd & $4$ & $4$\\ 
       & $0, 1, 2, 4, 8$ & $0, 0, 0, 1, 4$ & even & $5$ & $4$\\ 
       & $0, 2, 4, 6, 8$ & $0, 1, 2, 3, 4$ & odd & $10$ & $4$\\ \hline
       
 $6$ & $0, 1, 2, 3, 4, 10$ & $0, 0, 0, 0, 0, 5$ & odd & $5$ & $5$\\ 
       & $0, 1, 2, 3, 5, 10$ & $0, 0, 0, 0, 1, 5$ & even & $6$ & $5$\\ 
       & $0, 1, 2, 4, 6, 10$ & $0, 0, 0, 1, 2, 5$ & even & $8$ & $7$ \\       
       & $0, 1, 2, 5, 6, 10$ & $0, 0, 0, 2, 2, 5$ & odd & $9$ & $6$ \\              
       & $0, 2, 4, 6, 8, 10$ & $0, 1, 2, 3, 4, 5$ & odd & $15$ & $5$\\ \hline      
 \end{tabular}
 \end{center}

 \subsection{Genus 2}\label{genus2}
 Every genus $2$ curve is hyperelliptic, so every curve has exactly $2g+2=6$ subcanonical points, the ramification points of the hyperelliptic double cover.  These subcanonical points all have ramification sequence $0,1$, and the associated theta characteristics are odd.
 \subsection{Genus 3}\label{genus3}
 In genus $3$, we expect each component of $\scl_3$ to have codimension $3-1=2$ in $\mathcal{M}_{3,1}$, so each component of the locus of curves which have a subcanonical point should have codimension $1$ in $\mathcal{M}_{3}$.  This is of course true of the locus of hyperelliptic curves of genus $3$, which have subcanonical points with ramification sequence $0,1,2$.
 
 Any non-hyperelliptic curve of genus $3$ can be embedded by the canonical series as a smooth quartic in $\PR^2$, and conversely every smooth plane quartic is a non-hyperelliptic curve of genus $3$. 
 Given such a smooth plane quartic curve $C\subset \PR^2$, the canonical series is cut out on $C$ by the lines in $\PR^2$.  
 
 Thus a subcanonical point on $C$ is a point $p\in C$ such that for some line $L\subset \PR^2$ we have that $C\cap L = 4p$, or in other words $(C.L)_p=4$.  
  \begin{figure}[htb]
\centering 
\includegraphics[scale=.85]{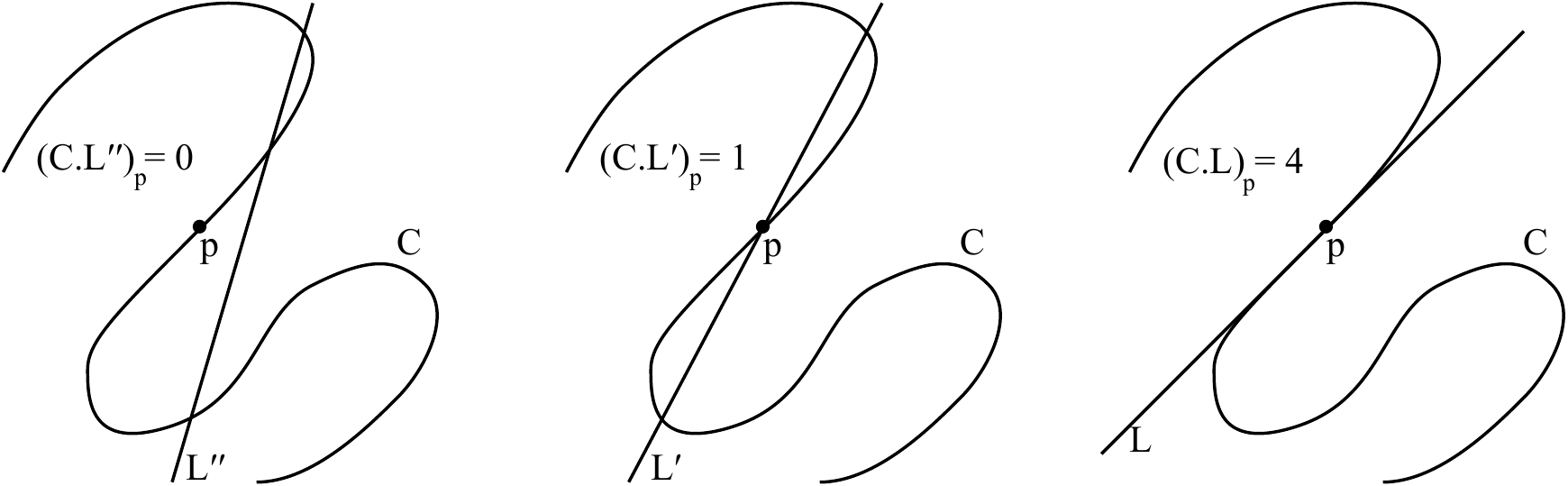} 
\caption{genus 3 non-hyperelliptic case} 
\end{figure}
In this case, any other line $L'\subset \PR^2$ either does not contain $p$, in which case $(C.L')_p=0$, or does contain $p$ but is not tangent to $C$ at $p$, so that $(C.L')_p=1$.  
 (If instead $(C.L')_p\geq 2$, then $(L.L')_p\geq \min \{(C.L')_p, (C.L)_p\}\geq 2$, and $L'=L$ by \bezout.)
 Thus the vanishing sequence for the canonical series of $C$ at $p$ is $0,1,4$, and the ramification sequence is $0,0,2$.

 The existence of such plane quartics is a simple application of Bertini's Theorem: if we fix a point and line $p\in L\subset\PR^2$, and consider the linear system of plane quartics $C$ satisfying the linear conditions $(C.L)_p\geq 4$, then this linear system has $p$ as its only basepoint (consider unions of four lines through $p$), so Bertini's Theorem shows that a general such $C$ is smooth away from $p$.  On the other hand, a general such $C$ is smooth at $p$ (for example, consider $C=L\cup C'$ with $p\not\in C'$), so we see that a general plane quartic satisfying $(C.L)_p\geq 4$ is smooth everywhere.  
 
 The locus of such plane quartics has dimension $\binom{6}{2}-1-4=10$, and the locus of choices of $p\in L$ has dimension $3$.  Two non-hyperelliptic curves of a given genus are isomorphic if and only if their canonical models are projectively equivalent.  Since a genus $3$ curve has only finitely many automorphisms and $\dim \PGL(3)=8$, we see that the locus in $\mathcal{M}_3$ of non-hyperelliptic genus $3$ curves which possess a subcanonical point has dimension $10+3-8=5$, or codimension $(3\cdot 3-3)-5=1$ in $\mathcal{M}_3$ as expected.  
 
 \par We might also have expected codimension one in this case since a general smooth plane quartic $C$ does have ordinary flexes, that is points $p$ whose tangent line $L$ satisfies $(C.L)_p=3$; possessing a hyperflex should be one additional condition.  Alternatively, note that a general smooth plane quartic has bitangents (in fact, it has $28$ bitangents, corresponding to the $28$ odd theta-characteristics), so it should be one additional condition for the two points of tangency to come together.   

 \subsection{Genus 4}\label{genus4}
 A non-hyperelliptic curve of genus $4$ canonically embeds as $C\subseteq\PR^3$, a degree $6$ curve which is the complete intersection of an irreducible quadric $Q$ and a cubic (cf. \cite{ACGH} ch. III).  
 
If $Q$ is smooth, then $C$ may be regarded as a smooth curve of bidegree $(3,3)$ on the surface $Q\isom \PR^1\times \PR^1$.  The canonical series $K_C$ is then cut on $C$ by hyperplanes in $\PR^3$, and thus by curves of bidegree $(1,1)$ on $Q$.  A curve of bidegree $(3,3)$ meets a curve of bidegree $(1,1)$ in $6$ points, counting multiplicities.  Suppose that $p$ is a subcanonical point of $C$, so that $C$ meets some bidegree $(1,1)$ curve $H$ at $p$ with multiplicity $6$ and nowhere else.  
\begin{figure}[htb]
\centering 
\includegraphics[scale=.51]{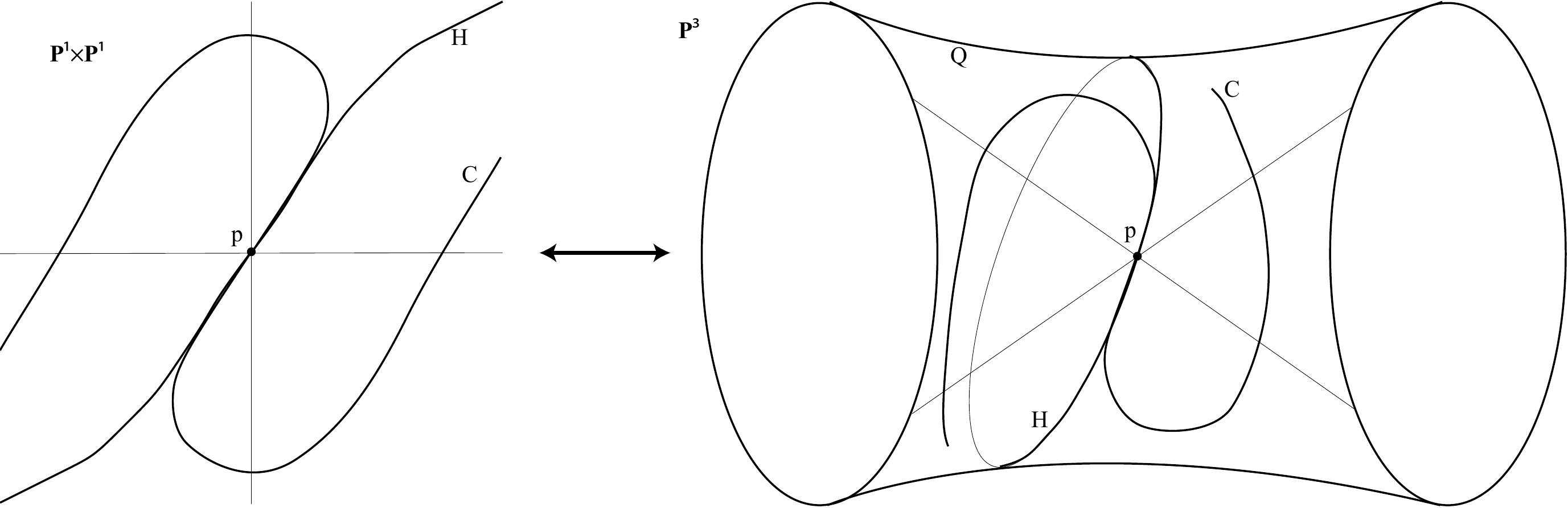} 
\caption{genus 4 odd case, vanishing sequence $0,1,2,6$} 
\end{figure}
We claim that if $H'$ is any other $(1,1)$ curve, then $(C.H')_p\leq 2$; for otherwise $(H.H')\geq 3$, but two $(1,1)$ curves meet in at most $2$ points unless they share a common component, and certainly $H$ must be irreducible, for if it consisted of two lines, they could not each meet $C$ at $p$ with multiplicity $3$.  We see then that the vanishing sequence of $K_C$ at $p$ must then be $0,1,2,6$, with corresponding ramification sequence $0,0,0,3$.  

To check that such curves do in fact exist, we note given a smooth quadric $Q\subseteq \PR^3$, so that $Q\isom \PR^1 \times \PR^1$, and a smooth curve $H\subset Q$ of bidegree $(1,1)$, we have that $H$ is a hyperplane section of $Q$, i.e. a smooth plane conic, so that $H\isom \PR^1$.  It may then be checked by direct computation that $\mathcal{O}_H(3)\isom \mathcal{O}_{\PR^1}(6)$ and that the restriction maps 
\[H^0(\PR^3,\mathcal{O}_{\PR^3}(3)) \rightarrow H^0(Q,\mathcal{O}_Q(3)) \isom H^0(\PR^1\times \PR^1, \mathcal{O}_{\PR^1\times \PR^1}(3,3)) \rightarrow H^0(\PR^1,\mathcal{O}_{\PR^1}(6))\]
are surjective.  Thus, given a point $p\in H$, there do exist curves of bidegree $(3,3)$ meeting $H$ at $p$ with multiplicity exactly $6$.  Moreover, meeting $H$ at $p$ with multiplicity at least $6$ imposes exactly $6$ linear conditions on the space of bidegree $(3,3)$ curves on $Q$.  By a Bertini's theorem argument, there exist smooth curves with this property.  By the above surjectivity, these curves arise as the complete intersection of a smooth quadric surface and a cubic surface, and one can check (e.g. using the adjunction formula) that such a curve does in fact have genus $4$.  We may estimate the dimension of the locus in $\mathcal{M}_{4,1}$ arising from these curves as follows:
\[\overbrace{\binom{5}{3}-1}^{\text{choice of $Q$}} + \overbrace{2\cdot 2-1}^{\text{choice of $H$}}+\overbrace{1}^{\text{$p$}} + \overbrace{4\cdot 4-1}^{\text{$C$ bidegree $(3,3)$}} - \overbrace{6}^{\text{$(C.H)_p\geq 6$}} - \overbrace{15}^{\text{$\PGL(4)$}}=7\]
 As expected, this estimate is equal to $\dim \scl_4^{\text{odd}}=2\cdot 4-1$.  A more rigorous version of this dimension count would show irreducibility as well. 
 
In order to find non-hyperelliptic curves of genus $4$ with even subcanonical points, we will need to look at the case where $Q$ is not smooth, but rather a cone over a smooth conic.  Here, we take a smooth curve $C$ on $Q$, defined by the intersection of $Q$ with a cubic hypersurface in $\PR^3$, so that it meets a point $p$ of a line $L$ of the ruling with multiplicity $3$.  Then the tangent plane $H$ to $Q$ at $p$ will intersect $Q$ in the double line $L$ and $H$ will intersect $C$ with multiplicity $6$, meaning $p$ is a subcanonical point of $C$.\footnote{To prove that such curves actually exist, we would first note that homogeneous cubic polynomials on $\PR^3$ restrict to $L$ to give all the homogeneous polynomials of degree $3$ on $L$, including those vanishing to order exactly $3$ at $p$.  We would then apply Bertini's theorem to the curves cut out on $Q$ by cubics that vanish to order at least $3$ along $L$.  Note that such curves containing the line $L$ certainly vanish to order at least $3$ along $L$ and can easily be made smooth at $p$.} 

\begin{figure}[htb]
\centering 
\includegraphics[scale=1.3]{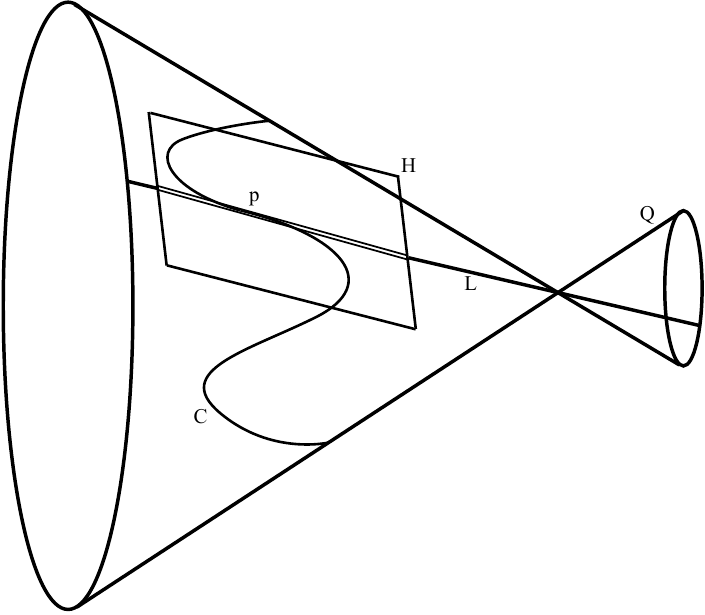} 
\caption{genus 4 even case, vanishing sequence $0,1,3,6$} 
\end{figure}
Here, however, if $H'$ is another hyperplane containing $L$, then $H'$ meets $C$ with multiplicity $3$.  We may certainly find hyperplanes not meeting $C$ at $p$ at all or meeting $C$ at $p$ with multiplicity $1$, so we see that the vanishing sequence for $K_C$ at $p$ is $0,1,3,6$, with corresponding ramification sequence $0,0,1,3$.
 
The dimension estimate of the corresponding locus in $\mathcal{M}_{4,1}$ is as follows:
\[\overbrace{\binom{5}{3}-1-1}^{\text{choice of cone $Q$}} + \overbrace{1}^{\text{$L$}}+\overbrace{1}^{\text{$p$}} + \overbrace{\binom{6}{3}-\binom{4}{3}-1}^{\text{cubic $C$ restricted to $Q$}} - \overbrace{3}^{\text{$(C.L)_p\geq 3$}} - \overbrace{15}^{\text{$\PGL(4)$}}=7.\]
As in the odd case, we see that this agrees with $\dim \scl_4^{\text{even}}= 2\cdot 4-1=7$.

 \subsection{Genus 5}\label{genus5}
 In genus $5$, the only possible ramification sequences for a non-hyperelliptic subcanonical point are still just $0,0,0,0,4$ and $0,0,0,1,4$; to show this, one can simply check that no other ramification sequences correspond to non-gap sequences that satisfy the semigroup condition.  We will describe a general genus $5$ curve with each of these ramification sequences.
 
First of all, any non-hyperelliptic, non-trigonal curve $\tilde{C}$ of genus $5$ is the normalization of a plane sextic $C$ with $5$ ordinary double points, $q_1,\ldots,q_5$.  The canonical series on $C$ is cut out by the plane cubics through the $5$ double points.  A point $p$ on $C$ is then subcanonical if there is some plane cubic $E$, passing through $q_1,\ldots, q_5$, which otherwise meets $C$ only at $p$, with multiplicity $8$.
\begin{figure}[htb]
\centering 
\includegraphics[scale=.65]{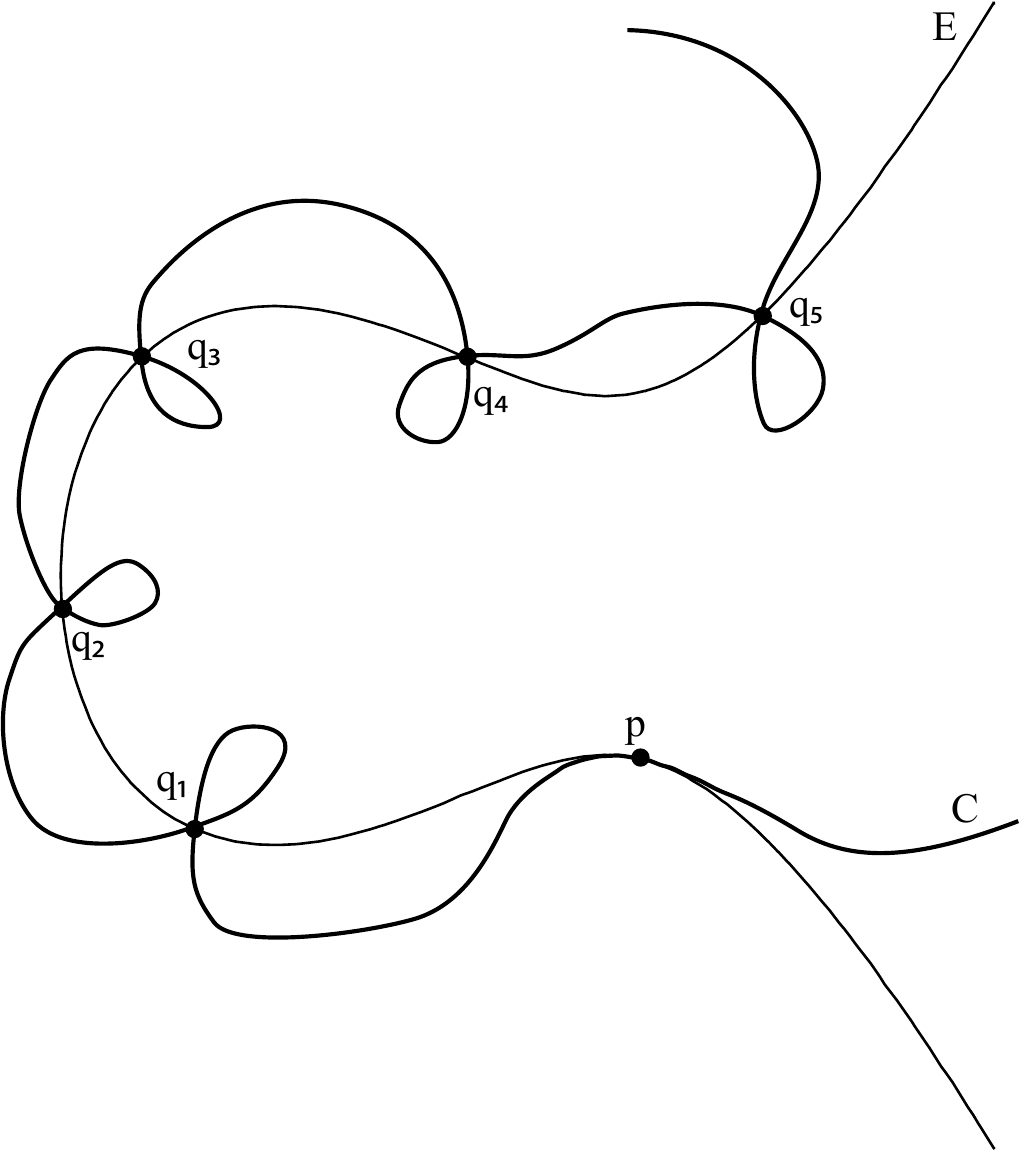} 
\caption{genus 5 non-hyperelliptic, non-trigonal case} 
\end{figure}
  Suppose that $E$ is a smooth elliptic curve.  Then the points $q_1,\ldots, q_5,p\in E$ are not arbitrary: since $C$ cuts out the divisor $2q_1+\ldots2q_5+8p$ on $E$, they must satisfy the relation \[2q_1+2q_2+\ldots+2q_5+8p\sim 6H\]
on $E$, where the divisor $H$ is cut out on $E$ by a line.

Now suppose that $E'\not=E$ is another cubic passing through $q_1,\ldots,q_5$ which satisfies $(E'.C)_p\geq 4$.  Then $(E.E')_p\geq 4$, so by \bezout, we must have $(E'.C)_p=(E.E')_p=4$, and in this case it must moreover be true that 
\[q_1+q_2+\ldots+q_5+4p\sim 3H\]
on $E$.  Conversely, if $q_1+q_2+\ldots+q_5+4p\sim 3H$, then there is a section of $\mathcal{O}_E(3)\isom \mathcal{O}_E(q_1+q_2+\ldots+q_5+4p)$ with zero divisor $q_1+q_2+\ldots+q_5+4p$, and since the map $H^0(\PR^2,\mathcal{O}_{\PR^2}(3))\to H^0(\PR^2,\mathcal{O}_{E}(3))$ is surjective (consider the short exact sequence 
$0\to \mathcal{I}_E(3)\to\mathcal{O}_{\PR^2}(3)\to\mathcal{O}_{E}(3)\to 0$) there is some cubic in $\PR^3$ which cuts out $q_1+q_2+\ldots+q_5+4p$ on $E$.

We thus see that if $q_1+q_2+\ldots+q_5+4p\sim 3H$ on $E$, then $p$ has vanishing sequence $0,1,2,4,8$ and ramification sequence $0,0,0,1,4$, whereas if $2q_1+2q_2+\ldots+2q_5+8p\sim 6H$ but $q_1+q_2+\ldots+q_5+4p\not\sim 3H$, then $p$ has vanishing sequence $0,1,2,3,8$ and ramification sequence $0,0,0,0,4$.

Of course, again we should check that such curves actually exist and calculate the dimensions of the corresponding loci in $\mathcal{M}_{5,1}$.  We begin by fixing a smooth cubic $E\subset \PR^2$ and points $q_1,\ldots,q_5,p\in E$ satisfying $2q_1+2q_2+\ldots+2q_5+8p\sim 6H$ on $E$.  We would like to compute how many linear conditions we impose on $H^0(\PR^2,\mathcal{O}_{\PR^2}(6))$ by requiring that a sextic curve $C$ have at least double points at $q_1,\ldots,q_5$ and satisfy $(C.E)_p\geq 8$.  We will then show by Bertini's theorem that there exist such curves which are smooth away from the $q_i$ and have simple nodes at the $q_i$.
 
  Now, since as above the map $H^0(\PR^2,\mathcal{O}_{\PR^2}(6))\to H^0(E,\mathcal{O}_{E}(6))$ is surjective, and since we have chosen $p$ and the $q_i$ so that $2q_1+2q_2+\ldots+2q_5+8p\sim 6H$, we find that there exist sextics $C\subset \PR^2$ which cut out exactly the divisor $2q_1+2q_2+\ldots+2q_5+8p$ on $E$.  Moreover,  
  we see that requiring that a sextic vanish along $E$ to order at least $2$ at each of the $q_i$ and to order $8$ at $p$ imposes exactly \[h^0(E,\mathcal{O}_{E}(6))-1=(18-1+1)-1=17\] linear conditions on $H^0(\PR^2,\mathcal{O}_{\PR^2}(6))$.  Now, for $C$ to have at least a double point at each point $q_i$ is $5$ additional linear conditions (that at each $q_i$ a derivative in a direction away from $E$ be zero as well); in fact, these $5$ linear conditions are independent, as we may consider sextics of the form $E\cup F_i$, where $F_i$ is a cubic vanishing at the four $q_j$ with $j\not= i$ but not at $q_i$.  
  
  Now, let $X$ be the projective space (of dimension $\binom{8}{2}-1-17-5$) of sextics which have double points (or worse) at the $q_i$ and which intersect $E$ at $p$ with multiplicity at least $8$.  A general $C\in X$ cuts out the divisor 
$2q_1+2q_2+\ldots+2q_5+8p$ on $E$.  We note first that the set of base points of $X$ is just $\{q_1,\ldots,q_5,p\}$, since for other points of $E$ we already know a general $C\in X$ does not contain them, and if $x\in \PR^2$ is not in $E$, we may find some plane cubic $E'$ which contains $q_1,\ldots,q_5$ but not $x$, and then $E\cup E'\in X$ would not contain $x$.  This shows, by Bertini's theorem, that a general $C\in X$ is smooth away from $\{q_1,\ldots,q_5,p\}$.

In fact, a general $C\in X$ must also be smooth at $p$, since if $C\not\supset E$ but $C$ is not smooth at $p$, then we may find some plane cubic $E'$ which contains $q_1,\ldots,q_5$ but not $p$, and then $C+EE'$ still satisfies the vanishing conditions and is smooth at $p$. (Here, in an abuse of notation, we are writing $C$, $E$, $E'$ for both the curves and their defining homogeneous polynomials.)  

Likewise, to control the singularity of $C$ at $q_i$, we may find a conic $Z$ containing $q_j$ for $j\not= i$ but not containing $q_i$, and then for a general choice of a line $L$ through $q_i$, the curves $C+EZL$ would still satisfy the vanishing conditions but have a simple node at $q_i$ whose branches have arbitrary tangent directions (aside from not being tangent to $E$, since then the curve would contain $E$).  This implies (e.g. by applying Bertini's theorem on the blowup of $\PR^2$ at the $q_i$) that a general $C\in X$ is smooth away from the $q_i$ with nodes at the $q_i$, as desired. 

To find the dimension of the corresponding loci (our parameter space will have two components depending on whether or not $q_1+q_2+\ldots+q_5+4p\sim 3H$), we must note that the map from an abstract smooth curve to the plane to give a degree $6$ singular curve is not unique, but rather there is a $2$-parameter family of such maps (essentially, the plane curves we are dealing with are projections of the canonical curve of degree $8$ in $\PR^4$ from two general points on the curve).  We thus calculate the dimensions of $\scl_5^{\text{odd}}$ and $\scl_5^{\text{even}}$ as 
\[\overbrace{\binom{5}{2}-1}^{\text{choice of $E$}} + \overbrace{6-1}^{2q_1+\ldots+2q_5+8p\sim 6H}+ \overbrace{\binom{8}{2}-1}^{\text{sextic $C$}} -\overbrace{17}^{\mathcal{O}_E(6)} - \overbrace{5}^{\text{$q_i$ double}} - \overbrace{8}^{\PGL(3)} -\overbrace{2}^{\dim G^2_6(C)}   =9.\]
This again agrees with the known dimension: $\dim \scl_5 = 2\cdot 5-1=9$.
     
 \subsection{Genus 6}\label{genus6} 
In genus $6$, there are more possible ramification sequences.  A curve $\tilde{C}$ of genus $6$ which is not hyperelliptic, trigonal, bielliptic, or isomorphic to a smooth plane quintic, is the normalization of a plane sextic $C$ having $4$ double points (cf. \cite{ACGH} V.A).  In this case, the situation is analogous to that in genus $5$: if $p\in C$ is the subcanonical point and $q_1,\ldots,q_4$ are the double points and $E\subset \PR^2$ is a smooth elliptic curve through the $q_i$ cutting out the divisor $10p$ on $C$, then \[2q_1+\ldots+2q_4+10p\sim 6H\] on $E$, where $H$ is a hyperplane section, and $p$ has vanishing sequence $0,1,2,3,5,10$ and ramification sequence $0,0,0,0,1,5$ if $q_1+\ldots+q_4+5p\sim 3H$ and otherwise has vanishing sequence $0,1,2,3,4,10$ and ramification sequence $0,0,0,0,0,5$.  

As in the genus $5$ case, we calculate the dimension as
\[\overbrace{\binom{5}{2}-1}^{\text{choice of $E$}} + \overbrace{5-1}^{2q_1+\ldots+2q_4+10p\sim 6H}+ \overbrace{\binom{8}{2}-1}^{\text{sextic $C$}} -\overbrace{17}^{\mathcal{O}_E(6)} - \overbrace{4}^{\text{$q_i$ double}} - \overbrace{8}^{\PGL(3)} -\overbrace{0}^{\dim G^2_6(C)}   =11\]
for the corresponding loci in $\mathcal{M}_{6,1}$, which again agrees with the known dimension of $2\cdot 6-1$ for $\scl_6^{\text{odd}}$ and $\scl_6^{\text{even}}$.

By Corollary \ref{bielliptic}, there exist subcanonical points on bielliptic curves of genus $6$ that have vanishing sequence $0,1,2,4,6,10$ and ramification sequence $0,0,0,1,2,5$.  To study the corresponding locus in $\mathcal{M}_{6,1}$, we recall from Section \ref{cyclic} that given an elliptic curve $E$ and distinct points $q_1,q_2,\ldots ,q_{10}\in E$, there exists a bielliptic double cover of $E$ of genus $6$ with branch locus $\{q_1,\ldots ,q_{10}\}$ having a subcanonical point that maps to $q_1$ if and only if $q_1+q_2+\ldots+q_{10}\sim 10q_1$ on $E$, and in this case the bielliptic double cover is unique.  We can show then that the corresponding locus in $\mathcal{M}_{6,1}$ is irreducible of dimension $\dim \mathcal{M}_{1,10}-1=9$.
 
A general subcanonical point $p$ on a smooth plane quintic $C$ has vanishing sequence $0,1,2,3,4,10$.  This is because the canonical series on $C$ is cut by plane conics, and if there is a smooth quadric $Q$ meeting $C$ only at $p$ with $(C.Q)_p=10$, then we can find quadrics meeting $C$ at $p$ with multiplicities $0,1,2,3,4$ simply by taking unions of lines (to get $(C,Q')_p=3$, for example, take the union of the tangent line to $C$ at $p$ with some other line through $p$).

 There however exist smooth plane quintics which possess a $5$-fold flex, that is, there exists a smooth plane quintic curve $C$ and a line $L\subset \PR^2$ meeting $C$ at a single point $p$ with $(C.L)_p=5$.  
   \begin{figure}[htb]
\centering 
\subfigure[$(C.Q)_p=2$]{\includegraphics[scale=.6]{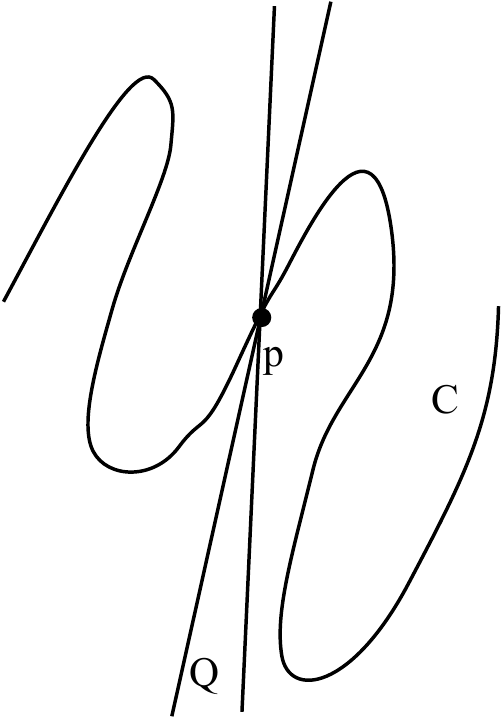}}\hspace{17pt}
\subfigure[$(C.Q)_p=5$]{\includegraphics[scale=.6]{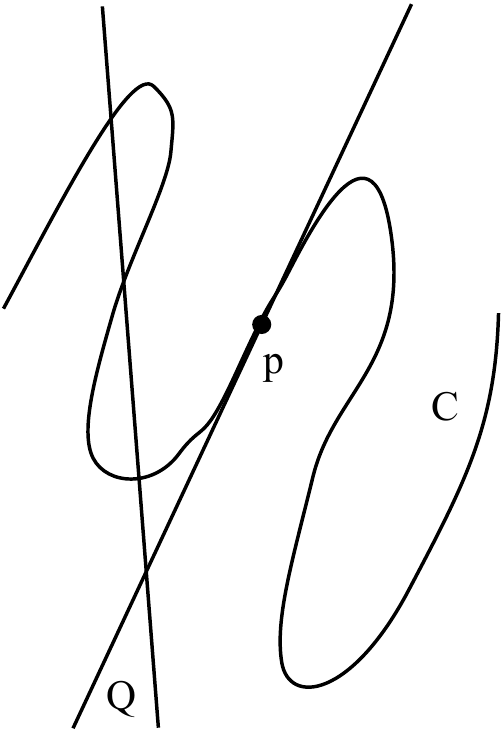}}\hspace{17pt}
\subfigure[$(C.Q)_p=6$]{\includegraphics[scale=.6]{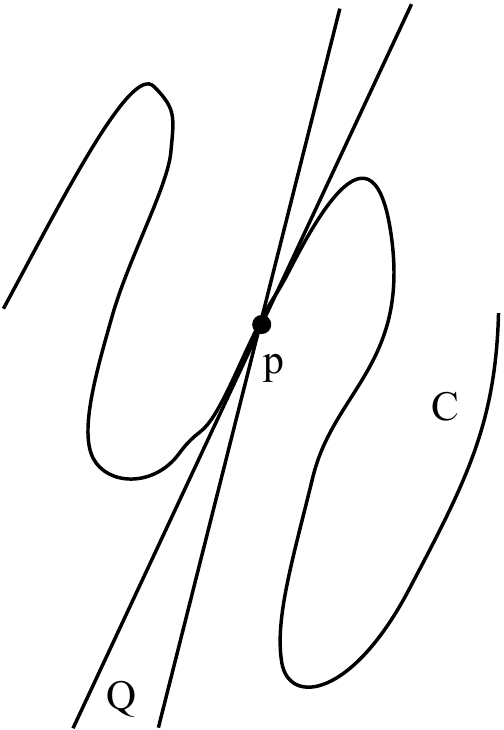}}\hspace{17pt}
\subfigure[$(C.Q)_p=10$]{\includegraphics[scale=.6]{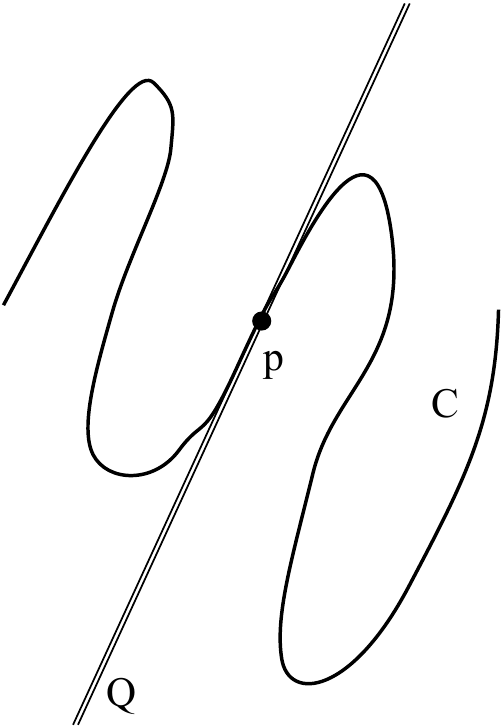}}
\caption{genus 6, plane quintic case} 
\end{figure}
 Then $p$ is a subcanonical point of $C$ with vanishing sequence $0,1,2,5,6,10$ and ramification sequence $0,0,0,2,2,5$.  To show this we can again construct sections with those vanishing orders by simply taking unions of lines.

 The proof that there exist smooth plane smooth quintics which possess a $5$-fold flex is entirely analogous to the proof in genus $3$ of the existence of smooth plane quartics with a hyperflex: we fix a point $p$ and a line $L$ containing it, and show by Bertini's theorem that there exists a smooth quintic $C$ satisfying $(C.L)_p=5$.  Since a smooth plane quintic can have only one embedding into $\PR^2$ (cf. \cite{ACGH} ch. V), we compute the dimension of the corresponding locus in $\mathcal{M}_{6,1}$ as:
\[\overbrace{2+1}^{\text{line and point}} + \overbrace{\binom{7}{2}-1}^{\text{quintic}}-\overbrace{5}^{\text{$5$-fold flex}} - \overbrace{8}^{\text{$\PGL(3)$}}-\overbrace{0}^{\dim G^2_5(C)} = 10.\]
A more detailed proof of this dimension count would also show the irreducibility of the corresponding locus.

\newpage
\singlespacing

\end{document}